\documentclass{article}
\usepackage{times}
\usepackage[margin=1.3in]{geometry}
\usepackage{helvet}  
\usepackage{courier}  
\usepackage{url}  
\usepackage{graphicx}  
\usepackage[utf8]{inputenc} 
\usepackage{booktabs}       
\usepackage{amsfonts}       
\usepackage{nicefrac}       
\usepackage{microtype}      
\usepackage{amsmath}
\usepackage{amssymb}
\usepackage{amsthm}
\usepackage{algorithm}
\usepackage{algorithmic}
\usepackage{subfigure}
\usepackage[usenames, dvipsnames]{color}
\usepackage{flexisym}

\DeclareMathOperator*{\argmin}{arg\ min}   
\DeclareMathOperator*{\prox}{prox}

\DeclareMathOperator*{\dom}{dom}
\DeclareMathOperator*{\with-probability}{with\ probability}

\DeclareMathOperator*{\Bernoulli}{Bernoulli}

\newcommand*{\QEDB}{\hfill\ensuremath{\square}}

\newtheorem{assumption}{\textbf{Assumption}}
\newtheorem{theorem}{\textbf{Theorem}}
\newtheorem{lemma}{\textbf{Lemma}}

\title{Double Quantization for Communication-Efficient Distributed Optimization}
\author{Yue Yu$^{1}$, Jiaxiang Wu$^{2}$, Longbo Huang$^{1}$\\
$^1$Institute for Interdisciplinary Information Sciences, Tsinghua University\\
$^2$Tencent AI Lab\\
yu-y14@mails.tsinghua.edu.cn,
jonathanwu@tencent.com,
longbohuang@tsinghua.edu.cn
}

\date{}

\begin{document}

\maketitle

\begin{abstract}
  Modern distributed training of machine learning models suffers from high communication overhead for synchronizing stochastic gradients and model parameters.
  In this paper, to reduce the communication complexity, we propose \emph{double quantization}, a general scheme for quantizing both model parameters and gradients.
  Three communication-efficient algorithms are proposed under this general scheme.
  Specifically, (i) we  propose a low-precision algorithm AsyLPG  with asynchronous parallelism,
  (ii) we explore integrating gradient sparsification with  double quantization and develop Sparse-AsyLPG,
  (iii)  we show that double quantization can also be accelerated by momentum technique and design accelerated AsyLPG.
  We establish rigorous performance guarantees for the algorithms, and
  conduct experiments on a multi-server test-bed to demonstrate that our algorithms can effectively save transmitted bits without performance degradation.
\end{abstract}

\section{Introduction}

%
Distributed optimization has received much recent attention due to data explosion and increasing  model complexity.
%
%
The \emph{data parallel} mechanism is a widely used distributed
architecture, which decomposes the time consuming gradient computations into sub-tasks, and assigns them to separate worker machines for execution.
%
Specifically, the training data is distributed to $M$ workers and each worker maintains a local copy of model parameters. 
At each iteration, each worker computes the gradient of a mini-batch randomly drawn from its local data. 
The global stochastic gradient is then computed by synchronously aggregating $M$ local gradients. 
Model parameters are then updated accordingly.

Two issues significantly slow down methods based on
 the data parallel architecture.
One is the communication complexity.
For example, all workers must send their entire local gradients to the master node at each iteration.
If gradients are dense, the master node has to receive and send $M*d$ floating-point numbers per iteration ($d$ is the size of the model vector), which scales linearly  with network  and model vector size.
With the increasing computing cluster size and model complexity, it has been  observed in many systems
that such communication overhead has become the performance bottleneck \cite{zhang2017zipml,wu2018error}.
The other factor is the synchronization cost, i.e.,
the master node has to wait for the last local gradient arrival  at each iteration.
This coordination dramatically increases system's idle time.

To overcome the first issue, many works focus on reducing the communication complexity of gradients in the data parallel architecture.
Generally, there are two kinds of methods.
One is quantization \cite{alistarh2017qsgd}, which stores gradients using fewer number of bits (lower precision).
%
The other is sparsification \cite{aji2017sparse}, i.e., dropping out some coordinates of gradients following certain rules.
%
%
However, existing communication-efficient algorithms based on data-parallel network still suffer from significant synchronization cost.
To address the second issue, many asynchronous algorithms have recently been developed for distributed training \cite{lian2015asynchronous,reddi2015variance}.
%
%
By allowing workers to communicate with master without synchronization, they can effectively improve the training efficiency.
Unfortunately, the communication bottleneck caused by transmitting gradients in floating-point numbers still exists.


In this paper, we are interested in jointly achieving  communication-efficiency and asynchronous parallelism. 
%
%
%
Specifically, we study the following composite problem
\begin{equation}
  \label{eq:pro-formu}
  \min_{x \in \Omega} P(x) = f(x) + h(x), \quad f(x) = \frac{1}{n}\sum\nolimits_{i=1}^n f_i(x),
\end{equation}
where $x \in \mathbb{R}^d$ is the model vector, $f_i(x)$ is smooth and $h(x)$ is convex but can be nonsmooth.
Domain $\Omega \subseteq \mathbb{R}^d$ is a convex set.
This formulation has found applications  in many different areas, such as multi-agent optimization \cite{nedic2009distributed} and distributed machine learning \cite{dean2012large}.

To solve (\ref{eq:pro-formu}), we propose 
a new \emph{double quantization} scheme,
which quantizes both model parameters and gradients. 
This quantization is nontrivial, because we have to deal with a low-precision gradient, which is evaluated on a low-precision model vector.
Three communication-efficient algorithms are then proposed under double quantization.
We analyze the precision loss of low-precision gradients and prove that these algorithms achieve fast convergence rate while significantly reducing the communication cost.
%
%
%
%
%
%
%
%
%
%
%
%
The main contributions are summarized as follows.

%
%
%

(i) We propose an \textbf{asy}nchronous \textbf{l}ow-\textbf{p}recision algorithm AsyLPG to solve the nonconvex and nonsmooth problem (\ref{eq:pro-formu}). 
We show that AsyLPG achieves the same asymptotic convergence rate as the unquantized serial algorithm, but with significantly less communication cost.

(ii)
We combine gradient sparsification with double quantization and propose
Sparse-AsyLPG to further reduce communication overhead. 
Our analysis shows that the convergence rate scales with $\sqrt{d/\varphi}$ for a sparsity budget $\varphi$.
%

%
%

(iii) 
We propose accelerated AsyLPG,
%
and mathematically prove that double quantization can be accelerated by momentum technique \cite{nesterov1983method,shang2017fast}.

(iv) We conduct experiments on a multi-server distributed test-bed. 
The results validate the efficiency of our algorithms. 
\begin{figure}
\begin{minipage}{0.5\textwidth}
  \centering
  \includegraphics[width=1.9in]{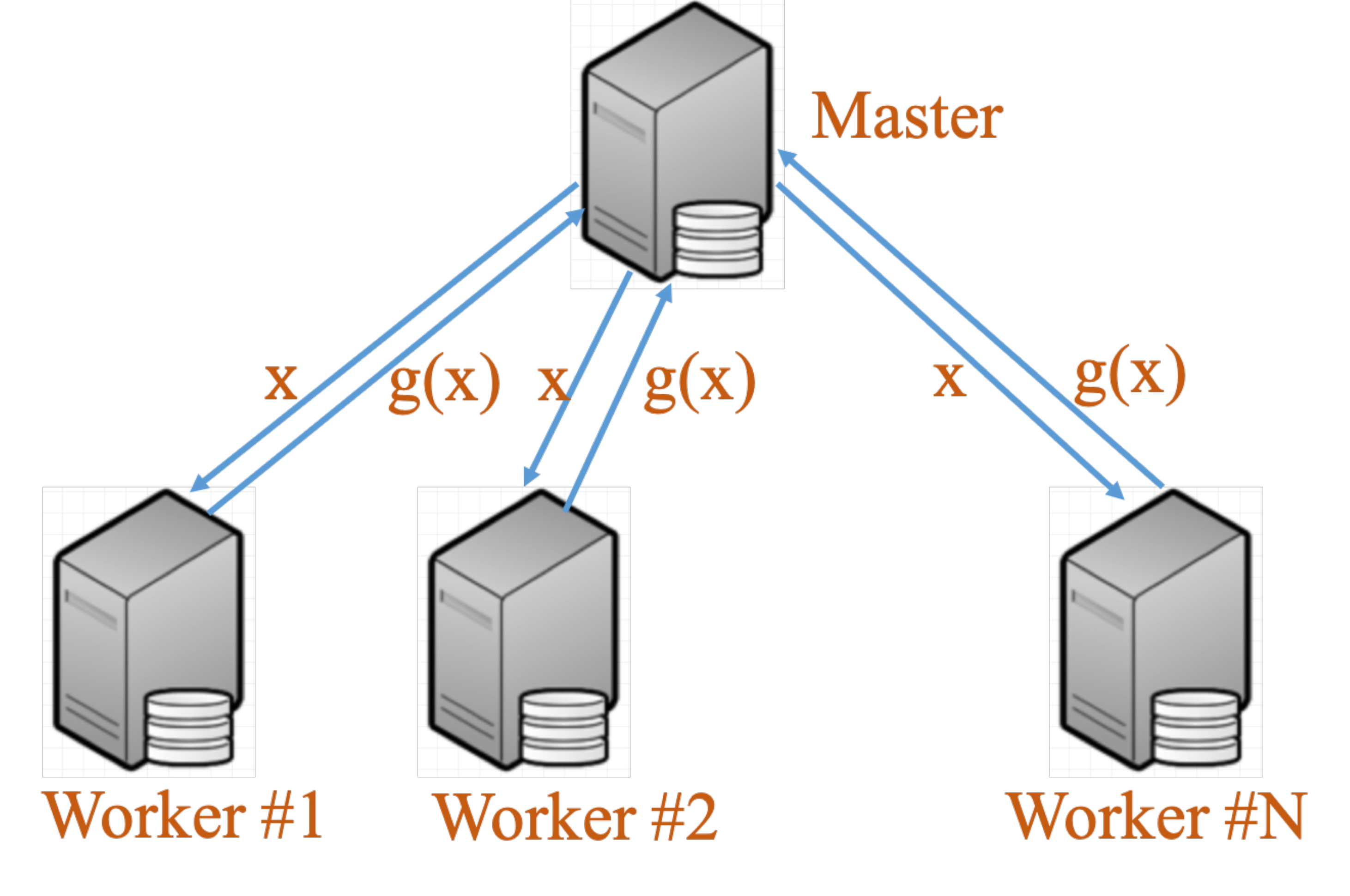}
\end{minipage}%
\begin{minipage}{0.5\textwidth}
  \centering
  \includegraphics[width=2.5in]{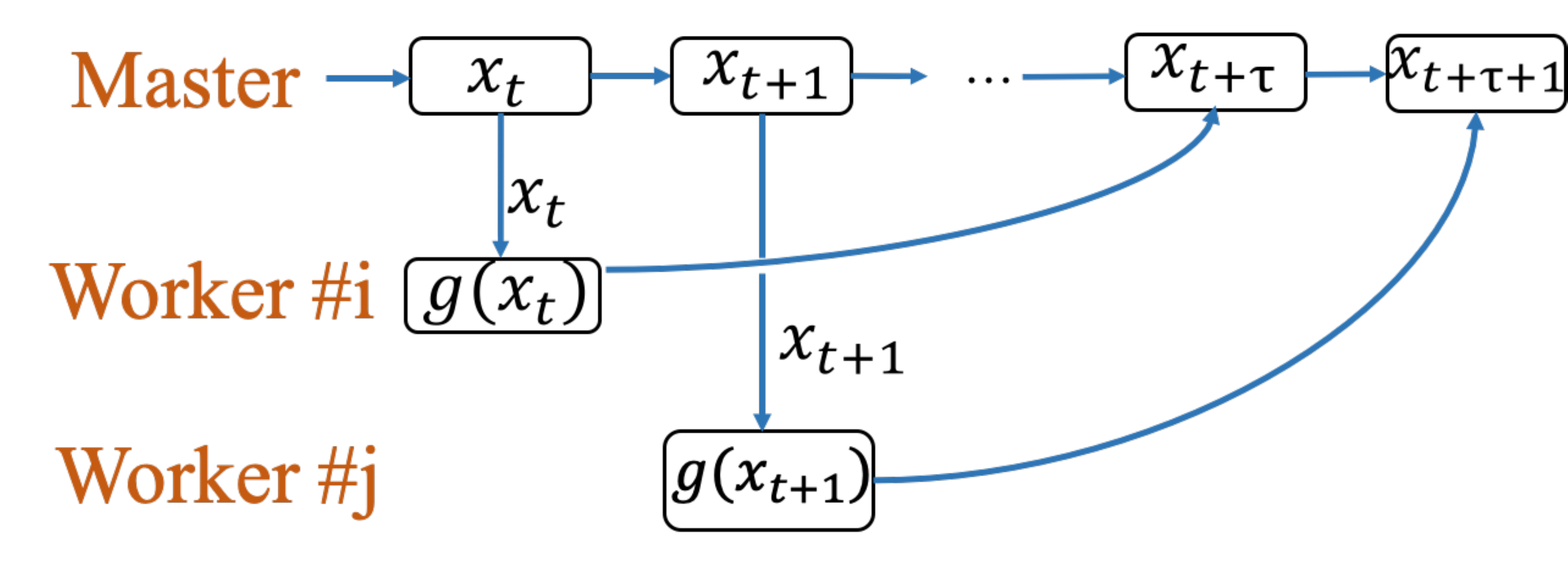}
\end{minipage}
\vskip -0.1in
\caption{The framework of distributed network with asynchronous communication. Left: network structure. Right: training process.}
\label{fig:asyn-network}
\vskip -0.1in
\end{figure}
\section{Related Work}
Designing large-scale distributed algorithms for machine learning has been receiving an increasing attention, and many algorithms, both synchronous and asynchronous, have been proposed, e.g., \cite{recht2011hogwild,bekkerman2011scaling,lian2015asynchronous,huo2016asynchronous}.
%
%
%
%
%
%
%
In order to reduce the communication complexity, researches also started to focus on
cutting down transmitted bits per iteration, based mainly on two schemes, i.e.,
%
quantization and sparsification.

\textbf{Quantization.}
Algorithms based on quantization  store a floating-point number using limited number of bits.
%
%
For example,
\cite{seide20141} quantized gradients to a representation of $\{-1,1\}$. 
\cite{alistarh2017qsgd,wen2017terngrad,wu2018error} adopted an unbiased gradient quantization with multiple levels.
The error-feedback method was applied in \cite{seide20141,wu2018error} to integrate history quantization error into the current stage.
%
%
\cite{de2018high} proposed a low-precision framework of SVRG \cite{johnson2013accelerating}, which quantized model parameters for single machine computation.
\cite{zhang2017zipml} proposed an end-to-end low-precision scheme, which quantized data, model and gradient with synchronous parallelism.
A biased quantization with gradient clipping was analyzed in \cite{yu2019AC}.
%
\cite{de2017understanding} empirically studied asynchronous and  low-precision SGD on logistic regression.

\textbf{Sparsification.}
Methods on sparsification drop out certain coordinates of gradients.  \cite{aji2017sparse} only transmitted gradients exceeding a threshold.
\cite{wangni2017gradient,wang2018atomo} formulated gradient sparsification into an optimization problem to balance sparsity and variance.
Recently, \cite{stich2018sparsified,alistarh2018convergence} analyzed the convergence behavior of sparsified SGD with memory, i.e., compensating gradient with sparsification error.
%
\vskip 0.02in
Our work distinguishes itself from the above results in:
%
(i) We quantize 
both
model vectors and gradients.
%
%
(ii) We integrate gradient sparsification within double quantization and prove the convergence.
(iii) We analyze how double quantization can be accelerated to reduce communication rounds.
%
%
\\
\\
\noindent
\textbf{Notation.}
$x^*$ is the optimal solution of (\ref{eq:pro-formu}).
%
 $||x||_{\infty}$, $||x||_1$ and $||x||$ denote the max, $L_1$ and $L_2$ norm of $x$.
For a vector $v_t \in \mathbb{R}^d$, $[v_t]_i$ or $v_{t,i}$ denotes its $i$-th coordinate.
$\{ \mathit{e}_i \}_{i=1}^d$ is the standard basis in $\mathbb{R}^d$.
%
The base of logarithmic function is $2$. 
%
$\tilde{O}(f)$ denotes $O(f\cdot polylog(f))$.
%
We use the proximal operator to handle the nonsmooth $h(x)$,
i.e., $\prox_{\eta h}(x) = \argmin_y h(y)+ \frac{1}{2\eta}||y-x||^2$.
If problem (\ref{eq:pro-formu}) is nonconvex and nonsmooth, 
we apply the commonly used convergence metric \emph{gradient mapping} \cite{nesterov2003introductory}: 
%
$G_{\eta}(x) \triangleq \frac{1}{\eta} [x - \prox_{\eta h}(x - \eta\nabla f(x) ) ]$.
%
%
$x$ is defined as an $\epsilon$-accurate solution if it satisfies $\mathbf{E}||G_{\eta}(x)||^2 \leq \epsilon$.
\section{Preliminary}
\subsection{Low-Precision Representation via Quantization}
\label{sec:quan}
%
Low-precision representation stores numbers using limited number of bits, contrast to the $32$-bit full-precision.\footnote{We assume without loss of generality that a floating-point number is stored using $32$ bits (also see \cite{alistarh2017qsgd,yu2019AC}). Our results can extend to the case when numbers are stored with other precision.}
It can be represented by a tuple $(\delta, b)$, where $\delta \in \mathbb{R}$ is the scaling factor and $b \in \mathbf{N^+}$ is the number of bits used.
Specifically, given a tuple $(\delta,b)$, the set of representable numbers is
\begin{equation*}
  \dom(\delta,b) = \{-2^{b-1}\cdot\delta,...,-\delta, 0, \delta,...,(2^{b-1}-1)\cdot\delta\}.
\end{equation*}
For any full-precision $x \in \mathbb{R}$, we call the procedure of transforming it to a low-precision representation as quantization, which is denoted by function $Q_{(\delta,b)}(x)$.
It outputs a number in $\dom(\delta,b)$ according to the following rules:

(i) If $x$ lies in the convex hull of $\dom(\delta,b)$, i.e., there exists a point $z \in \dom(\delta,b)$ such that $x \in [z, z+\delta]$,
then $x$ will be stochastically rounded in an unbiased way:
  \vskip -0.15in
\begin{equation*}
  Q_{(\delta,b)}(x) =
  \begin{cases}
    z+\delta, \quad    \with-probability  \ \frac{x-z}{\delta},\\
    z, \quad \quad \ \ \ \with-probability  \ \frac{z+\delta-x}{\delta}.
  \end{cases}
\end{equation*}
\vskip -0.05in
(ii) 
Otherwise, $Q_{(\delta,b)}(x)$  outputs the closest point to $x$ in $\dom(\delta,b)$, i.e., the minimum or maximum value.

This quantization method is widely used in existing works, e.g.,  \cite{zhang2017zipml,alistarh2017qsgd,wen2017terngrad,de2018high}, sometimes under different formulation.
In the following sections, we adopt $Q_{(\delta,b)}(v)$ to denote quantization on  vector $v \in \mathbb{R}^d$, 
which means that each coordinate of $v$ is independently quantized using the same tuple $(\delta,b)$.
Low-precision representation can effectively reduce communication cost, because we only need ($32+bd$) bits to transmit quantized  $Q_{(\delta,b)}(v)$ ($32$ bits for $\delta$, and $b$ bits for each coordinate),
%
%
whereas it needs $32d$ bits for a full-precision $v$.
%
\subsection{Distributed Network with Asynchronous Communication}
\label{sec:asyn-net}
%
%
%
%
%
%
As shown in Figure \ref{fig:asyn-network}, we consider a network with one master and multiple workers, e.g., the parameter-server setting.
The master maintains and updates a model vector $x$, and keeps a training clock.
Each worker gets access to the full datasets and keeps a disjoint partition of data.
%
%
In each communication round, a worker retrieves $x$ from the  master,  evaluates the gradient $g(x)$, and then sends it back to the master.
Since  workers asynchronously pull and push data during the training process, at a time $t$, the master may use a delayed gradient calculated on a previous $x_{D(t)}$, where $D(t) \leq t$.
%
%
Many works showed that a near linear speedup can be achieved if the delay is reasonably moderate \cite{lian2015asynchronous,reddi2015variance}.
%
%

%
%
\begin{figure}
  \scalebox{1.0}{
\begin{minipage}{0.5\linewidth}
  \centering
  \includegraphics[width=1.9in]{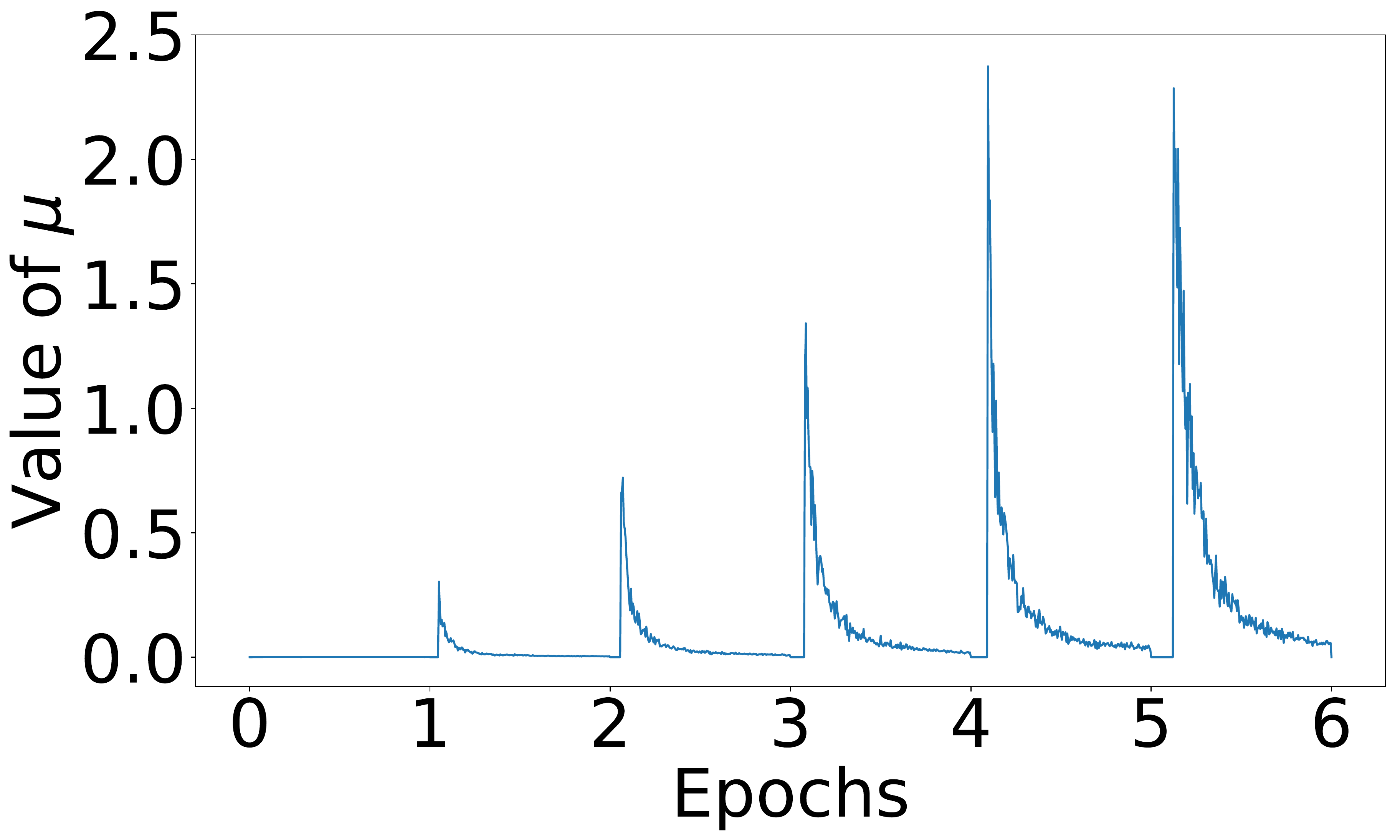}
\end{minipage}%
\begin{minipage}{0.5\linewidth}
  \centering
  \includegraphics[width=1.9in]{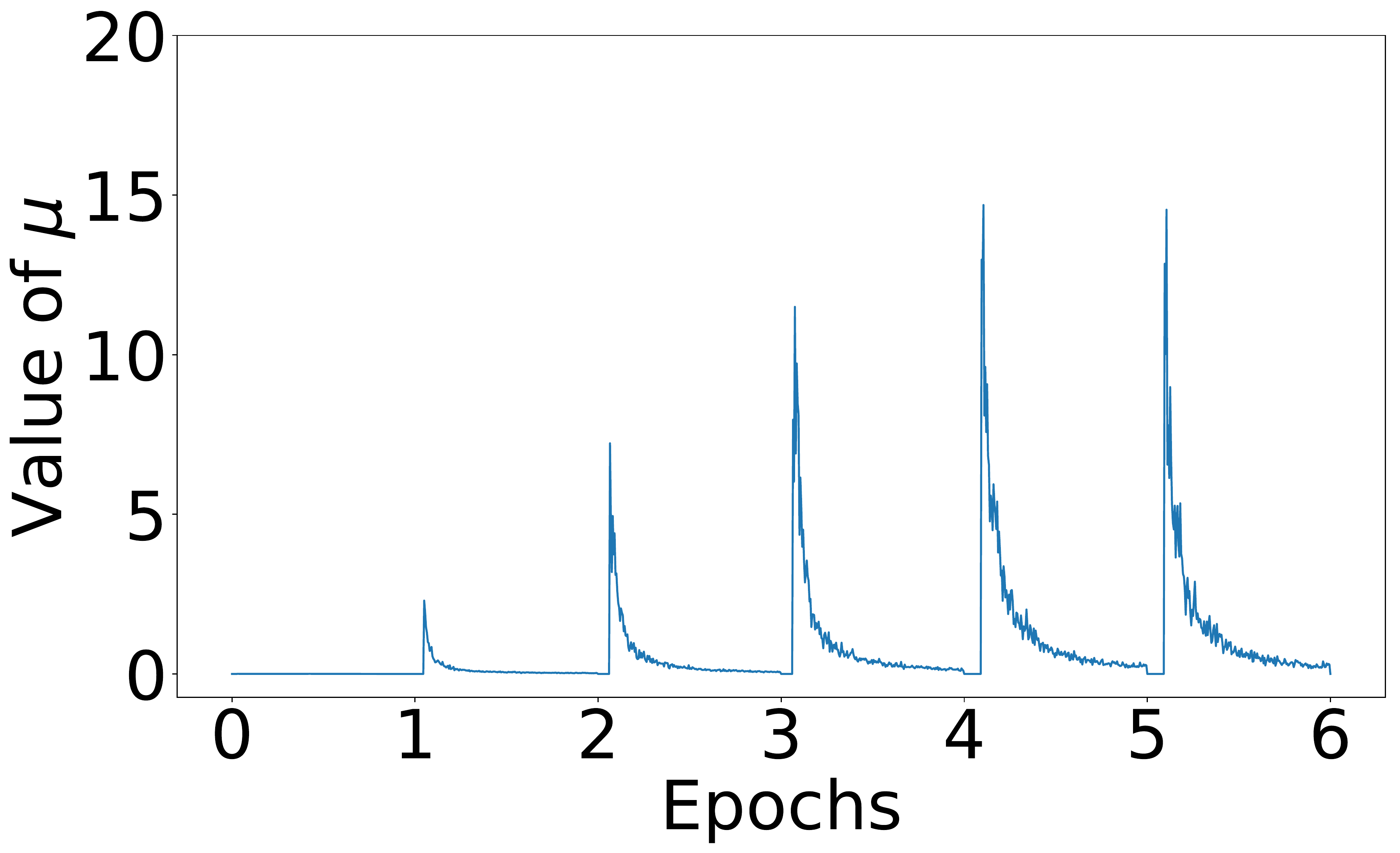}
\end{minipage}
}
\vskip -0.1in
\caption{The value of $\mu$ in different epochs that guarantees (\ref{eq:quan-x}). Left: $b_x = 8$. Right: $b_x = 4$. The statistics are based on a logistic regression on dataset \emph{covtype} ~\protect\cite{chang2011libsvm}.}
\label{fig:value-mu}
\vskip -0.1in
\end{figure}
\vskip -0.05in
\section{Algorithms}
\vskip -0.01in
To solve problem (\ref{eq:pro-formu}),
we propose a communication-efficient algorithm with double quantization, namely AsyLPG, and introduce its two variants with gradient sparsification and momentum acceleration.
%
%
We begin with the assumptions made in this paper. They are mild and are often assumed in the literature, e.g., \cite{johnson2013accelerating,Reddi2016Fast}.
\begin{assumption}
  \label{assump:unbiased}
The stochastically sampled gradient is unbiased, i.e., for
 $a\in\{1,...,n\}$ sampled in Step $11$ of Algorithms~\ref{Algm:AsyLPG}, \ref{Algm:Sparse-AsyLPG}, \ref{Algm:Acc-AsyLPG}, $\mathbf{E}_a[\nabla f_a(x)] = \nabla f(x)$.
Moreover, the random variables in different iterations are independent.
\end{assumption}
\begin{assumption}
  \label{assump:smooth}
  Each $f_i(x)$ in (\ref{eq:pro-formu}) is $L$-Lipschitz smooth, i.e., $||\nabla f_i(x) - \nabla f_i(y)|| \leq L ||x-y||$, $\forall x, y \in \Omega$.
\end{assumption}
\begin{assumption}
  \label{assump:delay}
  The gradient delay is upper bounded by some finite constant $\tau > 0$, i.e., $t- D(t) \leq \tau$, $\forall t$.
\end{assumption}
%
%
\begin{algorithm}[tb]
   \caption{AsyLPG}
\begin{algorithmic}[1]   \label{Algm:AsyLPG}
   \STATE {\bfseries Input:} $S$, $m$,  $\eta$, $b_x$, $b$, $\tilde{x}^0 = x^0$;
   \FOR{$s=0,1,...,S-1$}
       \STATE $x_0^{s+1} = \tilde{x}^s$;
       \STATE /* \textbf{Map-reduce global gradient computation}
       \STATE Compute $\nabla f(\tilde{x}^s) = \frac{1}{n}\sum\limits_{i=1}^n \nabla f_i(\tilde{x}^s)$;
       \FOR{$t=0$ {\bfseries to} $m-1$}
           \STATE /* \textbf{For master:}
           \STATE  (i) \textbf{Model Parameter Quantization:}
           Set $\delta_x = \frac{||x_{D(t)}^{s+1}||_{\infty}} {2^{b_x-1}-1}$ and quantize $x_{D(t)}^{s+1}$ subject to (\ref{eq:quan-x}). Then,
           send $Q_{(\delta_x, b_x)}(x_{D(t)}^{s+1})$ to workers;
           \STATE  (ii) Receive local gradient $Q_{(\delta_{\alpha_t}, b)}(\alpha_t)$,\\
           compute $u_t^{s+1} = Q_{(\delta_{\alpha_t}, b)}(\alpha_t) + \nabla f(\tilde{x}^s)$,  and \\
           update $x_{t+1}^{s+1} = \prox_{\eta h}(x_t^{s+1} - \eta u_t^{s+1})$;
           \medskip
           \STATE /* \textbf{For worker:}
           \STATE  (i) Receive $Q_{(\delta_x, b_x)}(x_{D(t)}^{s+1})$, stochastically sample a data-point $a \in \{1,...,n\}$, and calculate gradient
           $\alpha_t = \nabla f_a(Q_{(\delta_x, b_x)}(x_{D(t)}^{s+1})) - \nabla f_a(\tilde{x}^s)$;
           \STATE  (ii)
           \textbf{Gradient Quantization:} Set $\delta_{\alpha_t} = \frac{||\alpha_t||_{\infty}}{2^{b-1}-1}$ and send the quantized gradient $Q_{(\delta_{\alpha_t}, b)}(\alpha_t)$ to the master;
       \ENDFOR
      \STATE $\tilde{x}^{s+1} = x_m^{s+1}$;
  \ENDFOR
   \STATE {\bfseries Output:} Uniformly choosing from $\{\{x_t^{s+1}\}_{t=0}^{m-1}\}_{s=0}^{S-1}$.
\end{algorithmic}
\end{algorithm}
\subsection{Communication-Efficient Algorithm with Double Quantization: AsyLPG}
In this section, we introduce our new distributed algorithm AsyLPG, with \textbf{asy}nchronous communication and \textbf{l}ow-\textbf{p}recision floating-point representation.
As shown in Algorithm~\ref{Algm:AsyLPG},
AsyLPG divides the training procedure into epochs, similar to SVRG \cite{johnson2013accelerating}, with each epoch containing $m$ inner iterations.
%
%
%
%
At the beginning of each epoch, AsyLPG performs one round of communication between the master and workers to calculate the full-batch gradient $\nabla f(\tilde{x}^s)$, where $\tilde{x}^s$ is a snapshot variable evaluated at the end of each epoch $s$.
This one round communication involves full-precision operation because it is only performed once per epoch, and its communication overhead is small compared to the subsequent $m$ communication rounds in inner iterations, where the model parameters are updated.
%

In inner iterations, the communication between the master and workers utilizes the asynchronous parallelism described in Section~\ref{sec:asyn-net}.
To reduce communication complexity,   we propose \emph{double quantization}, i.e., quantizing both model parameters and gradients.

\textbf{Model Parameter Quantization.}
In Step $8$, before transmitting $x_{D(t)}^{s+1}$, the
master quantizes it into a low-precision format, subject to the following constraint:
\begin{equation}
  \label{eq:quan-x}
\mathbf{E}_Q||Q_{(\delta_x, b_x)}(x_{D(t)}^{s+1}) - x_{D(t)}^{s+1}||^2 \leq \mu||x_{D(t)}^{s+1} - \tilde{x}^s||^2.
\end{equation}
This condition is set to control the precision loss of
$x_{D(t)}^{s+1}$
with a positive hyperparameter $\mu$.
Note that with a larger $\mu$, we can aggressively save more transmitted bits (using a smaller $b_x$).
%
In practice, the precision loss and communication cost can be balanced by selecting a proper $\mu$ .
On the other hand, from the analysis aspect, we can always find a $\mu$ that guarantees (\ref{eq:quan-x}) throughout the training process, for any given $b_x$ (the number of quantization bits).
In the special case when $x_{D(t)}^{s+1} = \tilde{x}^s$, the master only needs to send a flag bit since $\tilde{x}^s$ has already been stored at the workers,
and (\ref{eq:quan-x}) still holds.
%
%
%
%
%
Figure~\ref{fig:value-mu} validates the practicability of (\ref{eq:quan-x}), where we plot the value of $\mu$ required to guarantee (\ref{eq:quan-x}). 
%
Note that the algorithm already converges in both graphs. In this case, we see that when $b_x$ is $4$ or $8$, $\mu$ can be upper bounded by a constant.
Also, by setting a larger $\mu$, we can choose a smaller $b_x$.
The reason $\mu$ increases at the beginning of each epoch is because $||x_{D(t)}^{s+1} - \tilde{x}^s||^2$ is small.
After several inner iterations, $x_{D(t)}^{s+1}$ moves further away from $\tilde{x}^s$. Thus, a smaller $\mu$ suffices to guarantee (\ref{eq:quan-x}).
%
%
%
%

\textbf{Gradient Quantization.}
After receiving the low-precision model parameter, as shown in Steps $11$-$12$, a worker calculates a gradient $\alpha_t$
and quantizes it into its low-precision representation $Q_{(\delta_{\alpha_t},b)} (\alpha_t)$, and then
sends it to the master.

In Step $9$, the master constructs a semi-stochastic gradient $u_t^{s+1}$ based on the received low-precision $Q_{(\delta_{\alpha_t},b)} (\alpha_t)$ and the full-batch gradient $\nabla f(\tilde{x}^s)$,
and updates the model vector $x$ using step size $\eta$.
The semi-stochastic gradient evaluated here adopts the variance reduction method proposed in SVRG \cite{johnson2013accelerating} and is used to accelerate convergence.
%
If Algorithm \ref{Algm:AsyLPG} is run without \emph{double quantization} and asynchronous parallelism, i.e., only one compute node with no delay, \cite{Reddi2016Fast} showed that:
\begin{theorem}
  \label{thm:reddi} (\cite{Reddi2016Fast}, Theorem $5$)
 Suppose $h(x)$ is convex and Assumptions \ref{assump:unbiased}, \ref{assump:smooth} hold.
  Let $T=Sm$ and $\eta = \rho/L$ where $\rho \in (0,1/2)$ and satisfies
  $4\rho^2 m^2 + \rho \leq 1$. Then for the output $x_{out}$ of Algorithm \ref{Algm:AsyLPG}, we have
  \begin{equation*}
    \mathbf{E}||G_{\eta}(x_{out})||^2 \leq \frac{2L(P(x^0) - P(x^*))}{\rho(1-2\rho)T}.
  \end{equation*}
\end{theorem}

\subsubsection{Theoretical Analysis}
\label{sec:AsyLPG-ana}
%
We begin with
the following lemma which bounds the variance of the delayed and low-precision gradient $u_t^{s+1}$.

\begin{lemma}
  \label{lemma:g-v}
  If Assumptions \ref{assump:unbiased}, \ref{assump:smooth}, \ref{assump:delay} hold,
  then for the gradient $u_t^{s+1}$ in Algorithm~\ref{Algm:AsyLPG}, its variance can be bounded by
  \begin{equation*}
    \begin{aligned}
    \mathbf{E}||u_t^{s+1} - \nabla f(x_t^{s+1})||^2 \leq
    2L^2(\mu+1)(\Delta+2)\mathbf{E}\Big[
    ||x_{D(t)}^{s+1} - x_t^{s+1}||^2 + ||x_t^{s+1} - \tilde{x}^s||^2 \Big],
\end{aligned}
  \end{equation*}
  where $\Delta = \frac{d}{4(2^{b-1}-1)^2}$.
\end{lemma}
\begin{theorem}
  \label{thm:asyLPG}
Suppose $h(x)$ is convex, conditions in Lemma~\ref{lemma:g-v} hold, $T=Sm$, $\eta = \frac{\rho}{L}$, where $\rho \in (0,\frac{1}{2})$, and $\rho, \tau$ satisfy
  \begin{equation}
    \label{eq:delay1}
     8\rho^2 m^2(\mu+1)(\Delta+2) + 2\rho^2(\mu+1)(\Delta+2)\tau^2 + \rho \leq 1.
  \end{equation}
Then, for the output $x_{out}$ of Algorithm \ref{Algm:AsyLPG} , we have
  \begin{equation*}
    \mathbf{E}||G_\eta(x_{out})||^2 \leq \frac{2L(P(x^0) - P(x^*))}{\rho(1-2\rho)T}.
  \end{equation*}
\end{theorem}
\textbf{Remarks.} From Theorem \ref{thm:asyLPG}, we see that if $\mu = O(1)$, $b = O(\log\sqrt{d})$ and $\rho = O(\frac{1}{m})$,
condition (\ref{eq:delay1}) can be easily satisfied, and AsyLPG
achieves the same asymptotic convergence rate as in Theorem~\ref{thm:reddi}, while transmitting much fewer bits ($b = O(\log{\sqrt{d}})$ is much smaller than $32$).
The comparisons between QSVRG \cite{alistarh2017qsgd} and AsyLPG in Figure~\ref{fig:real-sim}, \ref{fig:rcv1} and Table~\ref{fig:mnist-bits} also validate that our model parameter quantization significantly reduces transmission time and bits.
Note that AsyLPG can also adopt other gradient quantization methods (even biased ones), e.g., \cite{yu2019AC}, and similar results can be established.
\begin{algorithm}[tb]
   \caption{Sparse-AsyLPG}
\begin{algorithmic}[1]   \label{Algm:Sparse-AsyLPG}
   \STATE {\bfseries Input:} $S$, $m$, $b_x$, $b$, $\eta$, $\tilde{x}^0 = x^0$;
   \FOR{$s=0,1,...,S-1$}
       \STATE $x_0^{s+1} = \tilde{x}^s$;
       \STATE /* \textbf{Map-reduce global gradient computation}
       \STATE Compute $\nabla f(\tilde{x}^s) = \frac{1}{n}\sum\limits_{i=1}^n \nabla f_i(\tilde{x}^s)$;
       \FOR{$t=0$ {\bfseries to} $m-1$}
           \STATE /* \textbf{For master:}
           \STATE  (i) \textbf{Model Parameter Quantization:}
           Set $\delta_x = \frac{||x_{D(t)}^{s+1}||_{\infty}} {2^{b_x-1}-1}$ and quantize $x_{D(t)}^{s+1}$ subject to (\ref{eq:quan-x}).  Send $Q_{(\delta_x, b_x)}(x_{D(t)}^{s+1})$ to workers;
           \STATE  (ii) Receive local gradient $\zeta_t$,\\
           compute $u_t^{s+1} = \zeta_t + \nabla f(\tilde{x}^s)$ and \\
           update $x_{t+1}^{s+1} = \prox_{\eta h}(x_t^{s+1} - \eta u_t^{s+1})$;
          \medskip
           \STATE /* \textbf{For worker:}
           \STATE  (i) Receive $Q_{(\delta_x, b_x)}(x_{D(t)}^{s+1})$, stochastically sample a data-point $a \in \{1,...,n\}$, and calculate gradient
           $\alpha_t = \nabla f_a(Q_{(\delta_x, b_x)}(x_{D(t)}^{s+1})) - \nabla f_a(\tilde{x}^s)$;
          \STATE (ii) \textbf{Gradient Sparsification:} Choose a sparsity budget $\varphi_t$ and sparsify $\alpha_t$ to obtain $\beta_t$ using (\ref{eq:spar});
           \STATE  (iii)
           \textbf{Gradient Quantization:} Set $\delta_{\beta_t} = \frac{||\beta_t||_{\infty}}{2^{b-1}-1}$, and send  $\zeta_t= Q_{(\delta_{\beta_t}, b)}(\beta_t)$ to the master;
       \ENDFOR
      \STATE $\tilde{x}^{s+1} = x_m^{s+1}$;
  \ENDFOR
   \STATE {\bfseries Output:} Uniformly choosing from $\{\{x_t^{s+1}\}_{t=0}^{m-1}\}_{s=0}^{S-1}$.
\end{algorithmic}
\end{algorithm}
\subsection{AsyLPG with Gradient Sparsification}
%
%
%
In this section, we explore how to further reduce the communication cost
by incorporating gradient sparsification into double quantization.
As shown in Algorithm~\ref{Algm:Sparse-AsyLPG}, in Steps $11$-$13$, after calculating $\alpha_t$,  we successively perform sparsification and quantization on it.
Specifically, we drop out certain coordinates of $\alpha_t$ to obtain a sparsified vector $\beta_t$ according to the following rules \cite{wangni2017gradient}:
\vskip -0.05in
\begin{equation}
  \label{eq:spar}
  \beta_t = \Big[Z_1 \frac{\alpha_{t,1}}{p_1}, Z_2 \frac{\alpha_{t,2}}{p_2}, ..., Z_d \frac{\alpha_{t,d}}{p_d} \Big],
\end{equation}
%
where $Z = [Z_1, Z_2,..., Z_d]$ is a binary-valued vector with $Z_i \sim \Bernoulli(p_i)$, $0 < p_i \leq 1$, and $Z_i$'s are independent. 
Thus, $\beta_t$ is obtained by randomly selecting the $i$-th coordinate of $\alpha_t$ with probability $p_i$.
It can be verified that $\mathbf{E}[\beta_t] = \alpha_t$.
Define $\varphi_t \triangleq \sum_{i=1}^d p_i$ to measure the sparsity of $\beta_t$.
To reduce the communication complexity, it is desirable to make $\varphi_t$ as small as possible, which, on the other hand, brings about a large variance.
The following lemma quantifies the relationship between $\beta_t$ and $\varphi_t$, and is derived based on results in \cite{wang2018atomo}.

\begin{lemma}
  \label{lemma:spar-variance}
Suppose $\varphi_t \leq \frac{||\alpha_t||_1}{||\alpha_t||_{\infty}}$. Then, for $\alpha_t = \sum_{i=1}^d \alpha_{t,i} \mathit{e}_i$ and $\beta_t$ generated in (\ref{eq:spar}),
   we have $\mathbf{E}||\beta_t||^2 \geq \frac{1}{\varphi_t}||\alpha_t||^2_1$. 
   The equality holds if and only if $p_i = \frac{|\alpha_{t,i}| \cdot \varphi_t}{||\alpha_t||_1}$.
\end{lemma}
Based on Lemma~\ref{lemma:spar-variance}, in Step $12$ of Algorithm~\ref{Algm:Sparse-AsyLPG}, we can choose a sparsity budget $\varphi_t \leq \frac{||\alpha_t||_1}{||\alpha_t||_{\infty}}$ and set $p_i = \frac{|\alpha_{t,i}| \cdot \varphi_t}{||\alpha_t||_1}$ to minimize  the variance of $\beta_t$ (since $\beta_t$ is unbiased, its variance can be measured by the second moment).
%
%
%
After obtaining $\beta_t$, in Step $13$, we quantize it and send its low-precision version to the master.
The model parameter $x$ is updated in Step $9$, in a similar manner as in Algorithm~\ref{Algm:AsyLPG}.
In Algorithm~\ref{Algm:Sparse-AsyLPG}, we also employ model parameter quantization and asynchronous parallelism in inner iterations.

\subsubsection{Theoretical Analysis}
We first setup the variance of the sparsified gradient $u_t^{s+1}$.
\begin{lemma}
  \label{lemma:spar-var}
Suppose $\varphi_t \leq \frac{||\alpha_t||_1}{||\alpha_t||_{\infty}}$, Assumptions \ref{assump:unbiased}, \ref{assump:smooth}, \ref{assump:delay} hold, and for each $i \in \{1,...,d\}$, $p_i = \frac{|\alpha_{t,i}|\cdot \varphi_t}{||\alpha_t||_1}$.
  Denote $\Gamma = \frac{d^2}{4\varphi(2^{b-1}-1)^2} + \frac{d}{\varphi} +1$, where $\varphi = \min_t\{\varphi_t\}$.
  Then, for the gradient $u_{t}^{s+1}$ in Algorithm~\ref{Algm:Sparse-AsyLPG}, we have
  \begin{equation*}
    \begin{aligned}
    \mathbf{E}||u_t^{s+1} - \nabla f(x_t^{s+1})||^2 \leq
    2L^2(\mu+1)\Gamma \mathbf{E} \Big[ ||x_{D(t)}^{s+1} - x_t^{s+1}||^2 + ||x_t^{s+1} - \tilde{x}^s||^2 \Big].
\end{aligned}
  \end{equation*}
\end{lemma}
Based on the above variance bound, we state the convergence behavior of Sparse-AsyLPG in the following theorem.
\begin{theorem}
  \label{thm:sparse-asyLPG}
Suppose $h(x)$ is convex, conditions in Lemma~\ref{lemma:spar-var} hold, $T=Sm$, $\eta = \frac{\rho}{L}$, where $\rho \in (0,\frac{1}{2})$, and $\rho, \tau$ satisfy
  \begin{equation}
    \label{eq:delay3}
     8\rho^2 m^2(\mu+1)\Gamma + 2\rho^2(\mu+1)\tau^2\Gamma + \rho \leq 1.
  \end{equation}
Then, for the output $x_{out}$ of Algorithm \ref{Algm:Sparse-AsyLPG} , we have
  \begin{equation}
    \mathbf{E}||G_\eta(x_{out})||^2 \leq \frac{2L(P(x^0) - P(x^*))}{\rho(1-2\rho)T}.
  \end{equation}
\end{theorem}
\textbf{Remarks.} Setting $b = O(\log\sqrt{d})$, we obtain $\Gamma = O(d/\varphi)$.
We can conclude from Theorem~\ref{thm:sparse-asyLPG}  that Sparse-AsyLPG converges with a rate linearly scales with $\sqrt{d/\varphi}$,
and significantly reduces the number of transmitted bits per iteration.
Note that before transmitting $\zeta_t$, we need to encode it to a string, which contains $32$ bits for $\delta_{\beta_t}$ and $b$ bits for each coordinate.
Since $\beta_t$ is sparse, we only need to encode the nonzero coordinates, i.e.,
using $\log{d}$ bits to encode the position of a nonzero element  followed by its value.
\begin{algorithm}[tb]
   \caption{Acc-AsyLPG}
\begin{algorithmic}[1]   \label{Algm:Acc-AsyLPG}
   \STATE {\bfseries Input:} $S$, $m$, $b_x$, $b$, $\tilde{x}^0$, $y_m^0 = \tilde{x}^0$;
   \FOR{$s=1,2,...,S$}
       \STATE update $\theta_s, \eta_s$, $x_0^{s} = \theta_sy_0^s + (1-\theta_s)\tilde{x}^{s-1}$, $y_0^s = y_m^{s-1}$;
       \STATE /* \textbf{Map-reduce global gradient computation}
       \STATE Compute $\nabla f(\tilde{x}^{s-1}) = \frac{1}{n}\sum\limits_{i=1}^n \nabla f_i(\tilde{x}^{s-1})$;
       \FOR{$t=0$ {\bfseries to} $m-1$}
           \STATE /* \textbf{For master:}
           \STATE  (i) \textbf{Model Parameter Quantization:}
           Set $\delta_x = \frac{||x_{D(t)}^{s}||_{\infty}} {2^{b_x-1}-1}$, and quantize $x_{D(t)}^{s}$ subject to (\ref{eq:acc-quan-x}). Then,
           send $Q_{(\delta_x, b_x)}(x_{D(t)}^{s})$ to workers;
           \STATE  (ii) \textbf{Momentum Acceleration:}\\
           Receive local gradient $Q_{(\delta_{\alpha_t}, b)}(\alpha_t)$,\\
           compute $u_t^{s} = Q_{(\delta_{\alpha_t}, b)}(\alpha_t) + \nabla f(\tilde{x}^{s-1})$ and update\\
          $y_{t+1}^s = \prox_{\eta_s h}(y_t^s - \eta_s u_t^s)$,\\
           $x_{t+1}^{s} = \tilde{x}^{s-1} + \theta_s(y_{t+1}^s - \tilde{x}^{s-1})$;
           \medskip
           \STATE /* \textbf{For worker:}
           \STATE  (i) Receive $Q_{(\delta_x, b_x)}(x_{D(t)}^{s})$, stochastically sample a data-point $a \in \{1,...,n\}$ and calculate gradient
           $\alpha_t = \nabla f_a(Q_{(\delta_x, b_x)}(x_{D(t)}^{s})) - \nabla f_a(\tilde{x}^{s-1})$;
           \STATE  (ii)
           \textbf{Gradient Quantization:} Set $\delta_{\alpha_t} = \frac{||\alpha_t||_{\infty}}{2^{b-1}-1}$ and send the quantized gradient $Q_{(\delta_{\alpha_t}, b)}(\alpha_t)$ to the master;
       \ENDFOR
      \STATE $\tilde{x}^{s} = \frac{1}{m}\sum_{t=0}^{m-1}x_{t+1}^{s}$;
  \ENDFOR
   \STATE {\bfseries Output:} $\tilde{x}^S$.
\end{algorithmic}
\end{algorithm}
\subsection{Accelerated AsyLPG}
In the above, we mainly focus on reducing the communication cost within each iteration.
Here we propose
an algorithm
with an even faster convergence and
fewer
communication rounds.
Specifically, we incorporate the popular momentum or Nesterov technique \cite{nesterov1983method,shang2017fast} into AsyLPG. 
To simplify presentation, we only present accelerated AsyLPG (Acc-AsyLPG) in Algorithm~\ref{Algm:Acc-AsyLPG}. The method can similarly be applied to Sparse-AsyLPG.
%
%

 Algorithm~\ref{Algm:Acc-AsyLPG} still adopts asynchronous parallelism and  double quantization, and makes the following key modifications:
(i) in Step $8$, the model parameter quantization satisfy
\vskip -0.15in
\begin{equation}
  \label{eq:acc-quan-x}
    \mathbf{E}_Q||Q_{(\delta_x, b_x)}(x_{D(t)}^{s}) - x_{D(t)}^{s}||^2 \leq \theta_s\mu||x_{D(t)}^{s} - \tilde{x}^{s-1}||^2,
\end{equation}
%
where $\mu$ is the hyperparameter that controls the precision loss.
$\theta_s$ is the momentum weights and its value will be specified later.
(ii) Momentum acceleration is implemented in Steps $3$ and $9$, through an auxiliary variable $y_{t+1}^s$.
The update of $x_{t+1}^s$ combines history information $\tilde{x}^{s-1}$ and $y_{t+1}^s$.
%
%
%
%

In the following, we show  that with the above modifications, 
Acc-AsyLPG achieves an even faster convergence rate.
\begin{theorem}
  \label{thm:acc-AsyLPG}
Suppose each $f_i(x)$ and $h(x)$ are convex, Assumptions \ref{assump:unbiased}, \ref{assump:smooth}, \ref{assump:delay} hold, and the domain $\Omega$ of $x$ is bounded by D, such that $\forall x, y \in \Omega, ||x-y||^2 \leq D$.
  Let $\theta_s = \frac{2}{s+2}$, $\eta_s = \frac{1}{\sigma L \theta_s}$, where $\sigma > 1$ is a constant.
  %
  %
  If $\sigma$, $\tau$ satisfy
  $\tau \leq \frac{ \sqrt{ \big(\frac{2}{\gamma\theta_s} + \theta_s\Delta \big)^2 + \frac{4(\sigma -1)}{\gamma}} - (\frac{2}{\gamma\theta_s} + \theta_s\Delta)}{2}$
   where $\Delta = \frac{d}{(2^{b-1}-1)^2} + 2$, $\gamma = 1+2\theta_s\mu$, then under Algorithm~\ref{Algm:Acc-AsyLPG}, we have
  \begin{equation*}
    \mathbf{E}[P(\tilde{x}^S) - P(x^*)] \leq \tilde{O}((L/m + LD\mu\Delta/\tau + LD\mu)/S^2).
  \end{equation*}
\end{theorem}
The bounded domain condition in Theorem~\ref{thm:acc-AsyLPG} is commonly assumed in literature, e.g., \cite{wang2017memory}, and
the possibility of going outside domain is avoided by the proximal operator in Step $9$.
If $b = O(\log{\sqrt{d}})$ and $\mu = O(1)$, the constraint of delay $\tau$ can be easily satisfied with a moderate $\sigma$.
Then, our Acc-AsyLPG achieves acceleration while effectively reducing the communication cost.

\section{Experiments}
%
We conduct experiments to validate the efficiency of our algorithms.
%
The evaluations are setup on a $6$-server distributed test-bed. Each server has $16$ cores and $16$GB memory.
The communication between servers is handled by OpenMPI.\footnote{https://www.open-mpi.org/}
%
%

\begin{figure}
\begin{minipage}{0.5\textwidth}
  \centering
  \includegraphics[width=2.3in]{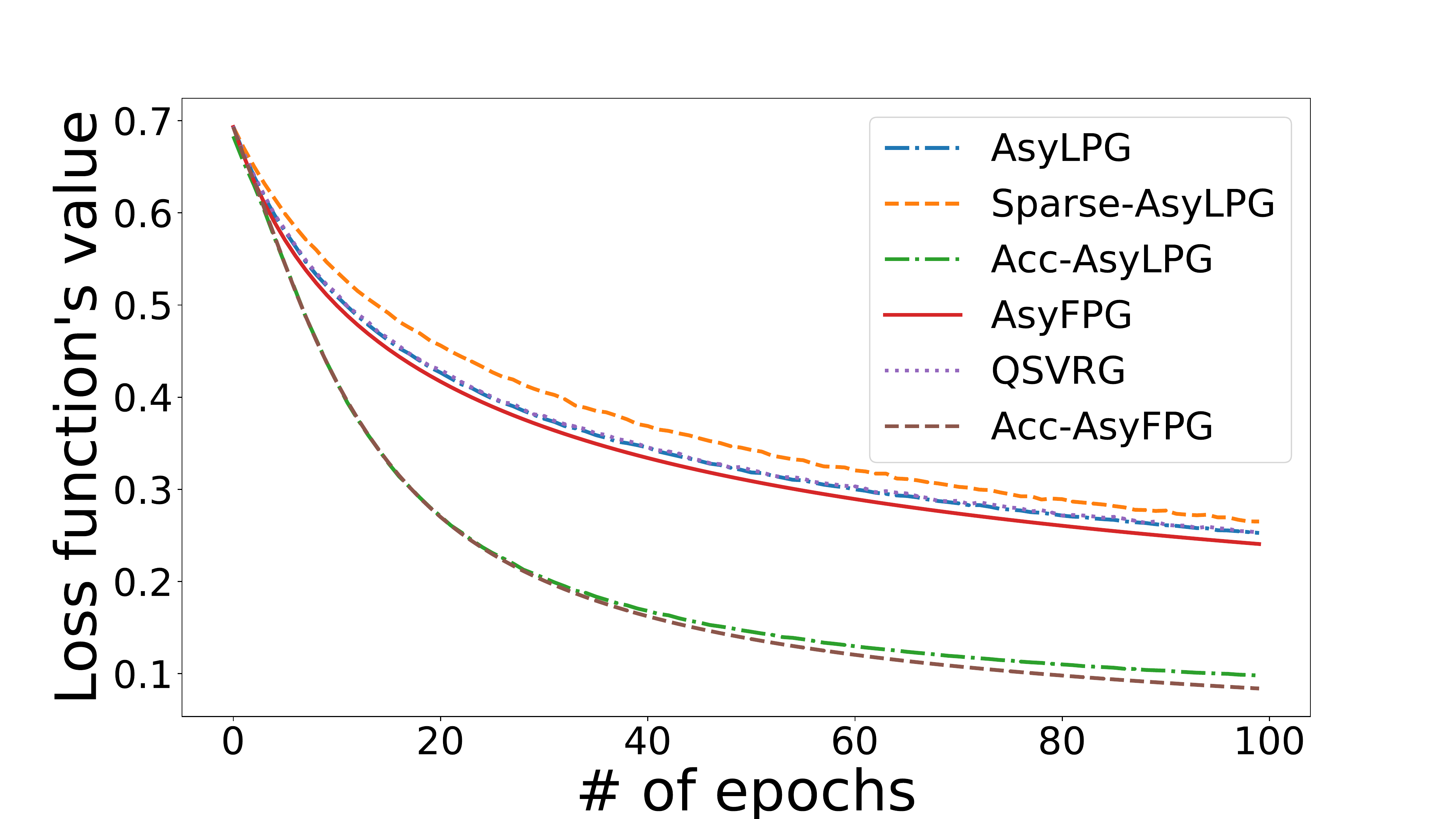}
\end{minipage}%
\begin{minipage}{0.5\textwidth}
  \centering
  \includegraphics[width=2.3in]{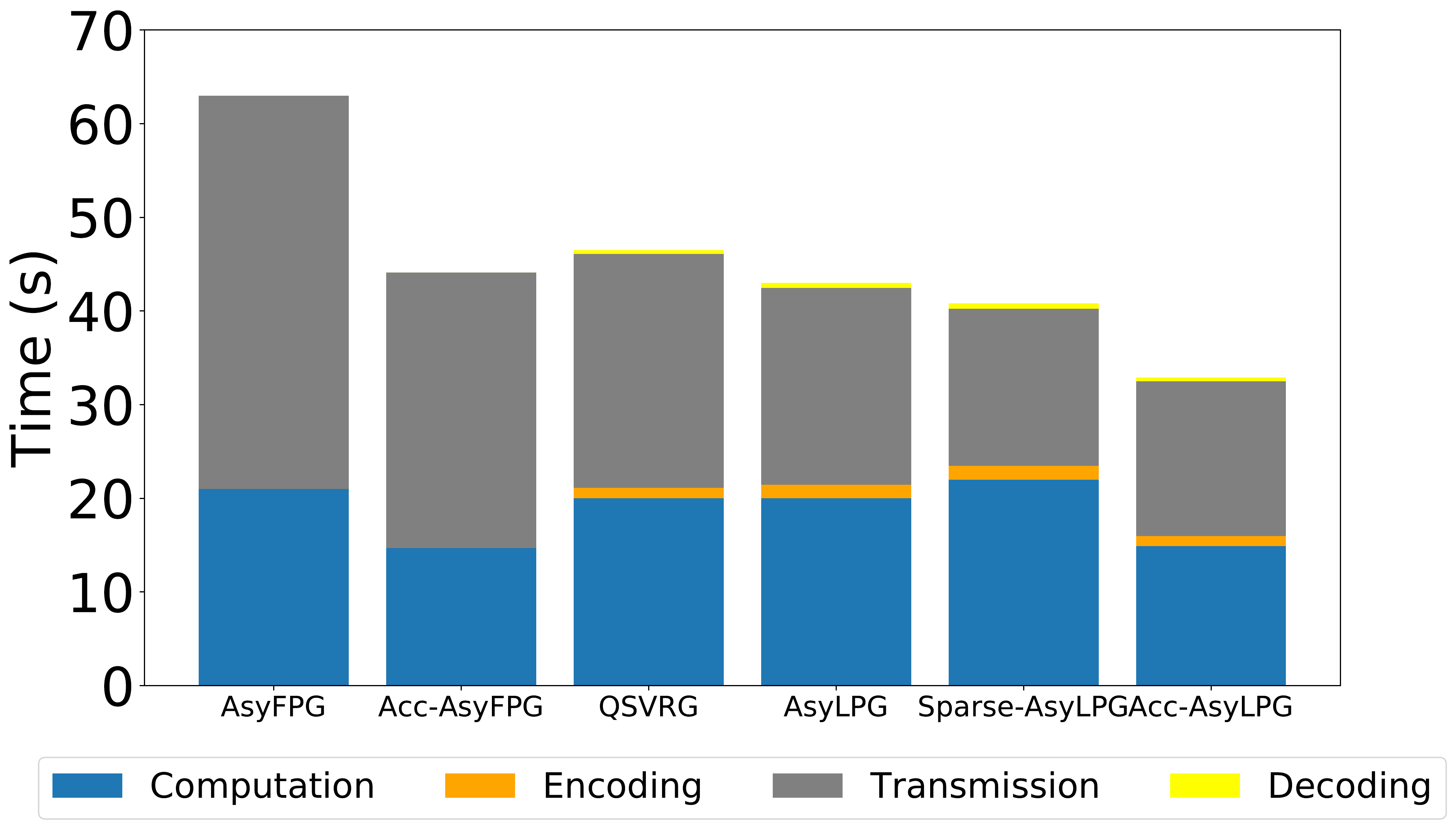}
\end{minipage}
\caption{Comparison of six algorithms on dataset \emph{real-sim}. Left: the training curve. Right: the decomposition of time consumption (the statistics are recorded until the training loss is first below $0.5$).}
\label{fig:real-sim}
\vskip -0.1in
\end{figure}
\begin{figure}
\begin{minipage}{0.5\linewidth}
  \centering
  \includegraphics[width=2.3in]{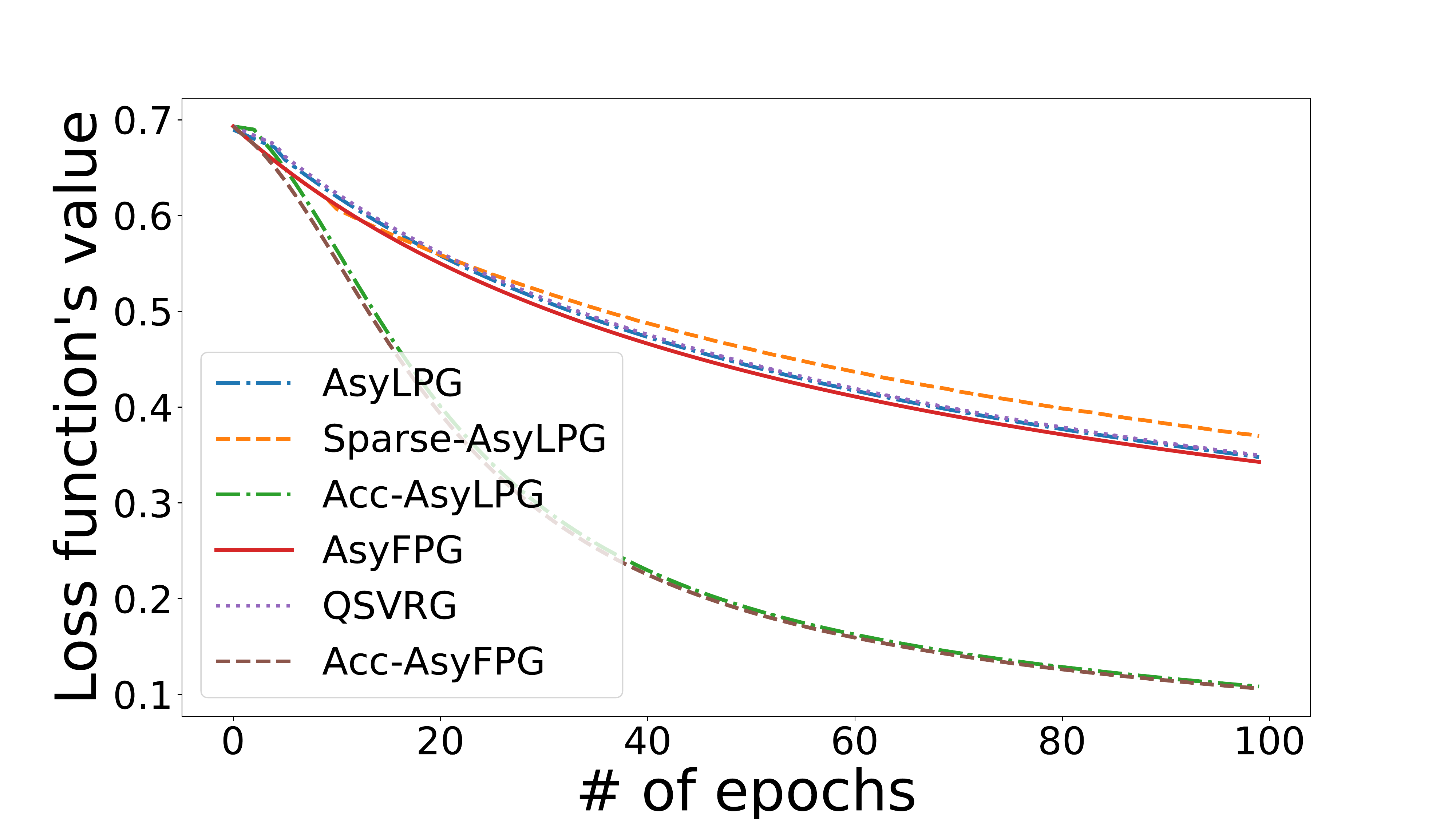}
\end{minipage}%
\begin{minipage}{0.5\linewidth}
  \centering
  \includegraphics[width=2.3in]{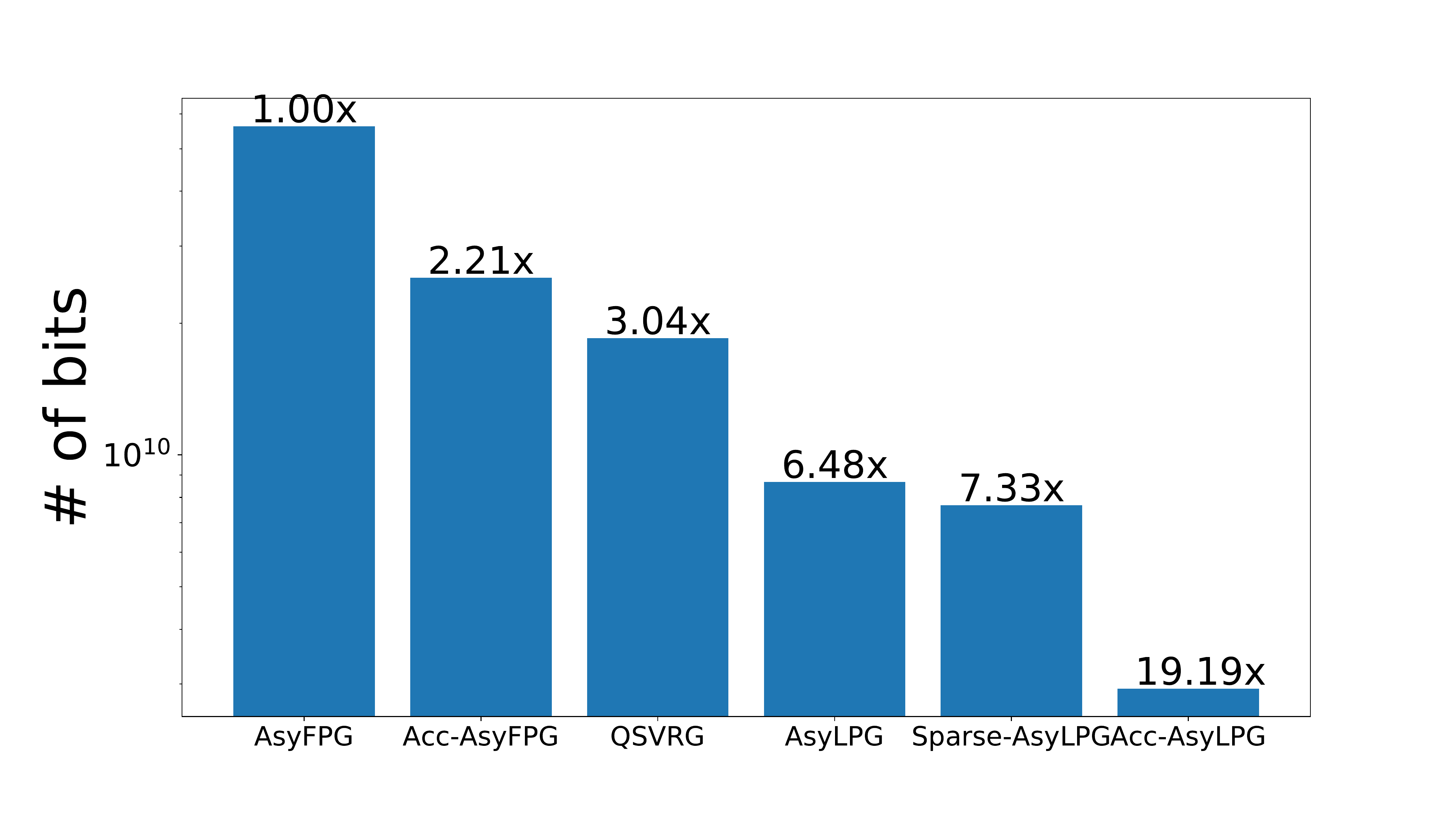}
\end{minipage}
\caption{Comparison on dataset \emph{rcv1}. Left: the training curve. Right: \# of transmitted bits until the training loss is first below $0.4$.}
\label{fig:rcv1}
\vskip -0.1in
\end{figure}
\subsection{Evaluations on Logistic Regression}
We begin with logistic regression on dataset \emph{real-sim} \cite{chang2011libsvm},
and use $L_1$, $L_2$ regularization with weights $10^{-5}$ and $10^{-4}$ respectively.
The mini-batch size $B = 200$ and epoch length $m = \lceil \frac{n}{B} \rceil$.
The following six algorithms are compared, using constant learning rate (denoted as \emph{lr})  tuned  to achieve the best result from $\{1e^{-1}, 1e^{-2}, 5e^{-2},1e^{-3}, 5e^{-3},...,1e^{-5},5e^{-5}\}$.

(i) AsyLPG, Sparse-AsyLPG, Acc-AsyLPG.
We set $b_x = 8$ and $b = 8$ in these three algorithms.
The sparsity budget in Sparse-AsyLPG is selected as $\varphi_t = ||\alpha_t||_1/||\alpha_t||_{\infty}$.
We do not tune $\varphi_t$ to present a fair comparison.
Parameters in Acc-AsyLPG are set to be $\theta_s = 2/(s+2)$ and $\eta_s = \emph{lr}/\theta_s$.

(ii) QSVRG \cite{alistarh2017qsgd}, which is a  gradient-quantized algorithm. 
We implement it in an asynchronous-parallelism way.
Its gradient quantization method is equivalent to Step $12$ in AsyLPG.
If run with synchronization and without quantization,
QSVRG and AsyLPG have the same convergence rate.
For a fair comparison, we set the gradient quantization bit $b = 8$ for QSVRG.

(iii) The full-precision implementations of AsyLPG and Acc-AsyLPG, denoted as AsyFPG and Acc-AsyFPG, respectively. In both algorithms, we remove double quantization.

\textbf{Convergence and time consumption.}
Figure~\ref{fig:real-sim} shows the evaluations on dataset \emph{real-sim}.
The left plot shows that AsyLPG and Acc-AsyLPG have similar convergence rates to their full-precision counterparts.
Our Sparse-AsyLPG also converges fast with little accuracy sacrifice.
%
%
%
%
The time consumption presented in the right plot shows the communication-efficiency of our algorithms.
%
%
With similar convergence rates, our low-precision algorithms significantly reduce the communication overhead when achieving the same training loss.
Moreover,
the comparison between AsyLPG and QSVRG
validates the redundancy of $32$ bits representation of model parameter.

\textbf{Communication complexity.}
We experimented logistic regression on dataset \emph{rcv1} \cite{chang2011libsvm}.
The $L_1$ and $L_2$ regularization are adopted, both with weights $10^{-4}$.
\emph{lr} is tuned in the same way as \emph{real-sim}.
In Figure~\ref{fig:rcv1}, we record the total number of transmitted bits of the six algorithms. 
It shows that AsyLPG, Sparse-AsyLPG, Acc-AsyLPG can save up to
6.48$\times$, 7.33$\times$ and 19.19$\times$ bits compared to AsyFPG.
%


\begin{figure}
\begin{minipage}{0.5\linewidth}
  \centering
  \includegraphics[width=2.3in]{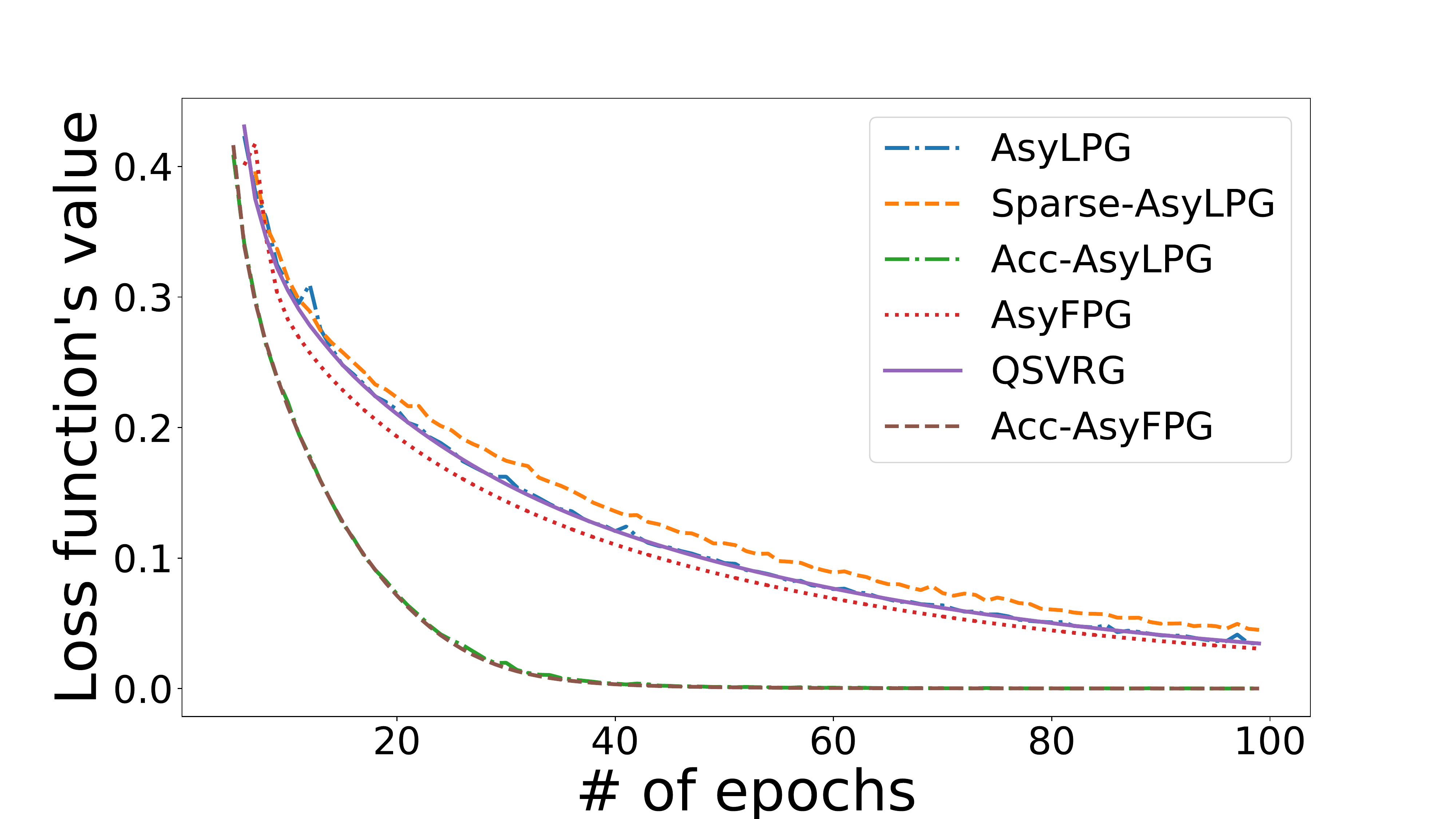}
\end{minipage}%
\begin{minipage}{0.5\linewidth}
  \centering
  \includegraphics[width=2.3in]{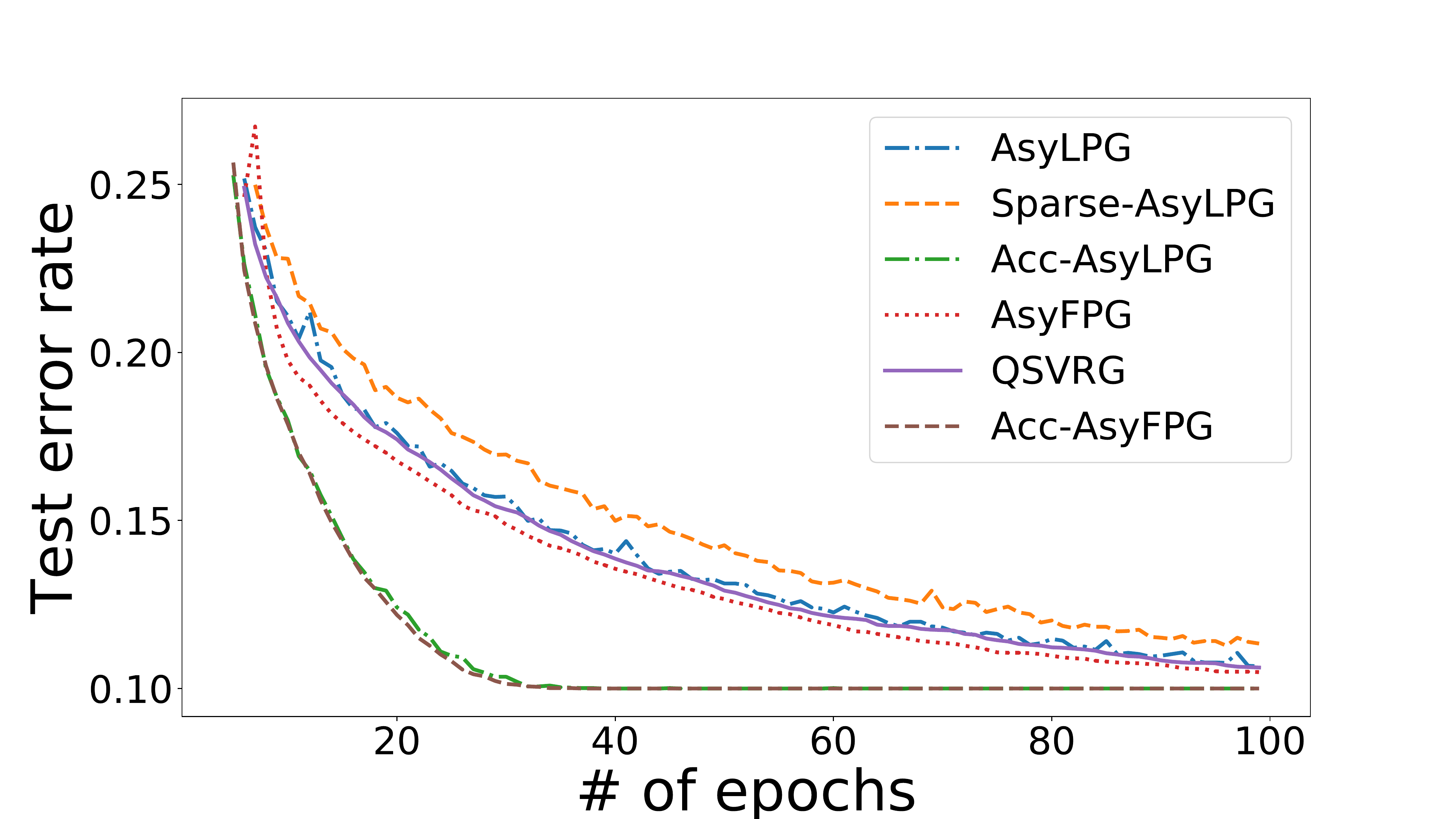}
\end{minipage}
\caption{Comparision of six algorithms on dataset MNIST. Left: the training curve, right: the test error rate.}
\label{fig:mnist-train}
\end{figure}

\begin{table}
\begin{minipage}{0.4\columnwidth}
  \centering
\scalebox{0.9}{
\begin{tabular}{lcc}
\toprule
Algorithm   & \# bits & Ratio \\
\midrule
AsyFPG      & 2.42$e$9  & $-$      \\
Acc-AsyFPG  & 6.87$e$8  & 3.52$\times$       \\
QSVRG            & 4.50$e$8  & 5.38$\times$      \\
AsyLPG    & 3.33$e$8  & 7.28$\times$      \\
Sparse-AsyLPG    & 2.73$e$8 & 8.87$\times$ \\
Acc-AsyLPG   & 1.26$e$8  & 19.13$\times$      \\
\bottomrule
\end{tabular}
}
\end{minipage}%
\begin{minipage}{0.6\columnwidth}
 \centering
  \scalebox{1.0}{
  \begin{tabular}{lcc| lcc}
  \toprule
  $\mu$   & $b_x$   &\#  bits &$\mu$   & $b_x$  & \# bits \\
  \midrule
  0.005    &   11  &   2.42$e$7   &  2.0    &  6 & 1.24$e$7     \\
  0.01    &   10   &  2.33$e$7    &  10    &  5 &1.14$e$7     \\
  0.05    &   9   & 1.52$e$7     & 50    &  4  &   1.08$e$7   \\
  0.1    &   8   &  1.07$e$7    &  150 & 3 & 1.44$e$7\\
  0.5    &   7  & 1.12$e$7   &  800 & 2 &  2.16$e$8 \\
  \bottomrule
  \end{tabular}
  }
\end{minipage}
\caption{Evaluation on dataset \protect\scalebox{0.9}{MNIST}. Left: \# of transmitted bits until the training loss is first below $0.05$. Right: The value of $b_x$ and \# of transmitted bits of AsyLPG under different $\mu$.}
\label{fig:mnist-bits}
\vskip -0.15in
\end{table}

%
\subsection{Evaluations on Neural Network}
We conduct evaluations on dataset MNIST\footnote{http://yann.lecun.com/exdb/mnist/}
using
a $3$-layer fully connected neural network.
%
The hidden layer contains $100$ nodes, and uses ReLU activation function.
Softmax loss function and $L_2$ regularizer with weight $10^{-4}$ are adopted. 
We use $10$k training samples and $2$k test samples
which are randomly drawn from the full dataset ($60$k training / $10$k test).
The mini-batch size is $20$ and the epoch size $m$ is $500$.
We set $b_x = 8$ and $b = 4$ for low-precision algorithms.
\emph{lr} is constant and is tuned to achieve the best result for each algorithm.

Figure~\ref{fig:mnist-train} presents the convergence  
of the six algorithms.
%
%
%
%
%
In the left table of Table \ref{fig:mnist-bits}, we record the total number of transmitted bits. 
We see that the results are similar as in the logistic regression case,
i.e.,
our new low-precision algorithms
can significantly reduce communication overhead 
compared to their full-precision counterparts and QSVRG.

\textbf{Study of $\mu$.}
The hyperparameter $\mu$ is set to control the precision loss incurred by model parameter quantization.
Before, we fix $b_x$ to compare our algorithms with other methods without model parameter quantization.
In the right table of Table~\ref{fig:mnist-bits}, we study the performance of AsyLPG under different $\mu$.
The value $b_x$ is searched by guaranteeing (\ref{eq:quan-x}).
We provide the overall numbers of transmitted bits under different $b_x$ until the training loss is first less than $0.5$.
The results validate that with the increasing $\mu$, we can choose a smaller $b_x$, to save more communication cost per iteration.
The total number of transmitted bits decreases until a threshold $\mu = 0.5$,  beyond which significant precision loss happens and we need more training iterations for achieving the same accuracy.

\section{Conclusion}
We propose three communication-efficient algorithms for distributed training with asynchronous parallelism.
The key idea is quantizing both model parameters and gradients,  called double quantization.
We analyze the variance of low-precision gradients and show that our algorithms achieve the same asymptotic convergence rate as the full-precision algorithms, while transmitting much fewer bits per iteration.
%
We also incorporate gradient sparsification into double quantization, and setup relation between convergence rate and sparsity budget.
We accelerate double quantization by integrating momentum technique.
The evaluations on logistic regression and neural network validate that our algorithms can significantly reduce communication cost.
\clearpage

%

%


\bibliography{sample_paper.bib}
\bibliographystyle{abbrv}

\clearpage

\section*{Supplementary Materials}
\section{Convergence Analysis for AsyLPG}
\begin{lemma}
  \label{lemma:unb-quan-var}
  For a vector $v \in \mathbb{R}^d$, if $\delta = \frac{||v||_{\infty}}{2^{b-1}-1}$ or $\frac{||v||_{2}}{2^{b-1}-1}$, we have
  \begin{equation}
    \mathbf{E} ||Q_{(\delta,b)}(v) - v||^2 \leq \frac{d\delta^2}{4}.
    \end{equation}
\end{lemma}
\noindent
\emph{Proof.} Because the squared $L_2$ norm separates along dimensions and each coordinate of $v$ is independently quantized, we only need to prove
$\mathbf{E}||Q_{(\delta,b)}([v]_i) - [v]_i||^2 \leq \frac{\delta^2}{4},$ for all $i \in \{1,...,d\}$.
If the scaling factor $\delta = \frac{||v||_{\infty}}{2^{b-1}-1}$ or $\frac{||v||_{2}}{2^{b-1}-1}$, it can be verified that $[v]_i$ locates in the convex hull of $\dom(\delta,b)$ and $Q_{(\delta,b)}(v)$ is an unbiased quantization.
Then $\mathbf{E}||Q_{(\delta,b)}([v]_i) - [v]_i||^2 \leq \frac{\delta^2}{4}$ according to Lemma $1$ in \cite{yu2019AC}.
\QEDB
\begin{lemma}
  \label{lemma:unbia-g-v}
  If Assumptions $1, 2$, $3$ hold,
  then for the gradient $u_t^{s+1}$ in Algorithm $1$,  its variance can be bounded by
  \begin{equation}
    \begin{aligned}
      \mathbf{E}||u_t^{s+1} - \nabla f(x_t^{s+1})||^2 \leq
      2L^2(\mu+1)(\Delta+2)\mathbf{E}\Big[
      ||x_{D(t)}^{s+1} - x_t^{s+1}||^2 + ||x_t^{s+1} - \tilde{x}^s||^2 \Big],
\end{aligned}
  \end{equation}
  where $\Delta = \frac{d}{4(2^{b-1}-1)^2}$.
\end{lemma}
\noindent
\emph{Proof.}
\begin{equation}
  \begin{aligned}
    \label{eq:k2}
    &\mathbf{E} ||u_t^{s+1} - \nabla f(x_t^{s+1})||^2 \\
    &= \mathbf{E}||Q_{(\delta_{\alpha_t},b)}(\alpha_t) + \nabla f(\tilde{x}^s) - \nabla f(x_t^{s+1})||^2 \\
    & = \mathbf{E}||Q_{(\delta_{\alpha_t},b)}(\alpha_t) - \alpha_t + \alpha_t + \nabla f(\tilde{x}^s) - \nabla f(x_t^{s+1})||^2 \\
    & = \mathbf{E}||Q_{(\delta_{\alpha_t},b)}(\alpha_t) - \alpha_t||^2 + \mathbf{E}||\alpha_t + \nabla f(\tilde{x}^s) - \nabla f(x_t^{s+1})||^2\\
    &\leq \Delta \underbrace{\mathbf{E} ||\alpha_t||^2}_{T_1} + \underbrace{\mathbf{E}||\alpha_t + \nabla f(\tilde{x}^s) - \nabla f(x_t^{s+1})||^2}_{T_2},
  \end{aligned}
\end{equation}
where the third equality holds because $\delta_{\alpha_t} = \frac{||\alpha_t||_{\infty}}{2^{b-1}-1}$ and $Q_{(\delta_{\alpha_t},b)}(\alpha_t)$ is an unbiased quantization.
The final inequality follows from Lemma \ref{lemma:unb-quan-var}.
Next we bound $T_1$ and $T_2$.
\begin{equation}
  \begin{aligned}
    T_1 &= \mathbf{E}||\nabla f_a(Q_{(\delta_x, b_x)}(x_{D(t)}^{s+1})) - \nabla f_a(\tilde{x}^s)||^2\\
    &\leq L^2 \mathbf{E}||Q_{(\delta_x, b_x)}(x_{D(t)}^{s+1}) - \tilde{x}^s||^2\\
    & = L^2\mathbf{E}||Q_{(\delta_x, b_x)}(x_{D(t)}^{s+1}) - x_{D(t)}^{s+1}||^2 + L^2 \mathbf{E}||x_{D(t)}^{s+1} - \tilde{x}^s||^2\\
    & \leq L^2(\mu+1)\mathbf{E}||x_{D(t)}^{s+1} - \tilde{x}^s||^2 \\
    & \leq 2L^2(\mu+1) \mathbf{E}||x_{D(t)}^{s+1} - x_t^{s+1}||^2 + 2L^2(\mu+1)\mathbf{E}||x_t^{s+1} - \tilde{x}^s||^2,
  \end{aligned}
\end{equation}
where the first inequality adopts the Lipschitz smooth property of $f_a(x)$,
and the second equality holds because $\delta_x = \frac{||x_{D(t)}^{s+1}||_{\infty}} {2^{b_x-1}-1}$ and $Q_{(\delta_{x},b_x)}(x_{D(t)}^{s+1})$ is an unbiased quantization.
The second inequality uses the condition in Step $8$ of Algorithm $1$.
%
%
With similar arguments, we obtain the upper bound of $T_2$ in the following.

\begin{equation}
  \begin{aligned}
    T_2 &= \mathbf{E}||\nabla f_a(Q_{(\delta_x, b_x)}(x_{D(t)}^{s+1})) - \nabla f_a(\tilde{x}^s) + \nabla f(\tilde{x}^s) - \nabla f(x_t^{s+1})||^2\\
    &\leq 2\mathbf{E}||\nabla f_a(Q_{(\delta_x, b_x)}(x_{D(t)}^{s+1})) - \nabla f_a(x_{D(t)}^{s+1})||^2 + 2\mathbf{E}||\nabla f_a(x_{D(t)}^{s+1}) - \nabla f_a(\tilde{x}^s) + \nabla f(\tilde{x}^s) - \nabla f(x_t^{s+1})||^2\\
    &\leq 2\mathbf{E}||\nabla f_a(Q_{(\delta_x, b_x)}(x_{D(t)}^{s+1})) - \nabla f_a(x_{D(t)}^{s+1})||^2 + 4\mathbf{E}||\nabla f_a(x_{D(t)}^{s+1}) - \nabla f_a(x_t^{s+1})||^2 \\
     &\quad + 4\mathbf{E}||\nabla f_a(x_t^{s+1}) - \nabla f_a(\tilde{x}^s) + \nabla f(\tilde{x}^s) - \nabla f(x_t^{s+1})||^2\\
    &\leq 2\mathbf{E}||\nabla f_a(Q_{(\delta_x, b_x)}(x_{D(t)}^{s+1})) - \nabla f_a(x_{D(t)}^{s+1})||^2 + 4\mathbf{E}||\nabla f_a(x_{D(t)}^{s+1}) - \nabla f_a(x_t^{s+1})||^2 + 4\mathbf{E}||\nabla f_a(x_t^{s+1}) - \nabla f_a(\tilde{x}^s)||^2\\
    & \leq 2L^2\mathbf{E}||Q_{(\delta_x, b_x)}(x_{D(t)}^{s+1}) - x_{D(t)}^{s+1}||^2 + 4L^2\mathbf{E}||x_{D(t)}^{s+1} - x_t^{s+1}||^2 + 4L^2\mathbf{E}||x_t^{s+1} - \tilde{x}^s||^2\\
    &\leq 4L^2(\mu+1)\mathbf{E}||x_{D(t)}^{s+1} - x_t^{s+1}||^2 + 4L^2(\mu+1)\mathbf{E}||x_t^{s+1} - \tilde{x}^s||^2.
  \end{aligned}
\end{equation}
where in the third inequality we adopt $\mathbf{E}||\nabla f_a(x_t^{s+1}) - \nabla f_a(\tilde{x}^s) + \nabla f(\tilde{x}^s) - \nabla f(x_t^{s+1})||^2 \leq \mathbf{E}||\nabla f_a(x_t^{s+1}) - \nabla f_a(\tilde{x}^s)||^2$. It is true because $\mathbf{E}||x - \mathbf{E}[x]||^2 \leq \mathbf{E}||x||^2$.
The last inequality follows from Step $8$ of Algorithm $1$.
Putting them together, we obtain Lemma \ref{lemma:unbia-g-v}.
\QEDB
\\
\\
\noindent
\textbf{\emph{Proof of Theorem $2$.}}
Define $\bar{x}_{t+1}^{s+1} \triangleq \prox_{\eta h}(x_t^{s+1} - \eta \nabla f(x_t^{s+1}))$.
According to equations $(8)$-$(12)$ in \cite{Reddi2016Fast}, we get
\begin{equation}
  \begin{aligned}
  \label{eq:P-in}
  \mathbf{E} & \Big[P(x_{t+1}^{s+1})\Big]\\
  & \leq \mathbf{E} \Big[ P(x_t^{s+1})  + (L-\frac{1}{2\eta})||\bar{x}_{t+1}^{s+1} - x_t^{s+1}||^2 + (\frac{L}{2} - \frac{1}{2\eta})||x_{t+1}^{s+1} - x_t^{s+1}||^2 - \frac{1}{2\eta}||x_{t+1}^{s+1} - \bar{x}_{t+1}^{s+1}||^2 \\
  & \quad + \langle x_{t+1}^{s+1} - \bar{x}_{t+1}^{s+1}, \nabla f(x_t^{s+1}) - u_t^{s+1} \rangle \Big]\\
  & \leq \mathbf{E} \Big[ P(x_t^{s+1}) + \frac{\eta}{2}||u_t^{s+1} - \nabla f(x_t^{s+1})||^2 + (L-\frac{1}{2\eta})||\bar{x}_{t+1}^{s+1} - x_t^{s+1}||^2 + (\frac{L}{2} - \frac{1}{2\eta})||x_{t+1}^{s+1} - x_t^{s+1}||^2\Big].
\end{aligned}
\end{equation}
Using Lemma \ref{lemma:unbia-g-v}, we have
\begin{equation}
  \begin{aligned}
  \mathbf{E} \Big[P(x_{t+1}^{s+1}) \Big]
  \leq & \mathbf{E} \Big[ P(x_t^{s+1}) + \eta L^2(\mu+1)(\Delta + 2)||x_{D(t)}^{s+1} - x_t^{s+1}||^2 + \eta L^2(\mu+1)(\Delta + 2) ||x_t^{s+1} - \tilde{x}^s||^2 \\
  &+ (L-\frac{1}{2\eta})||\bar{x}_{t+1}^{s+1} - x_t^{s+1}||^2 + (\frac{L}{2} - \frac{1}{2\eta})||x_{t+1}^{s+1} - x_t^{s+1}||^2\Big].
\end{aligned}
\end{equation}
Define $R_t^{s+1} \triangleq \mathbf{E}\Big[P(x_t^{s+1}) + c_t||x_t^{s+1} - \tilde{x}^s||^2\Big]$,
where $\{c_t\}_{t=0}^{m}$ is a nonnegative decreasing sequence with $c_m = 0$,
 $c_t = c_{t+1}(1+\beta) + \eta L^2(\mu+1)(\Delta + 2)$ and $\beta = \frac{1}{m}$.
 Therefore,
 \begin{equation}
   \begin{aligned}
     \label{eq:bound-c}
     c_0 &\leq \eta L^2(\mu+1)(\Delta + 2) \cdot \frac{(1+\beta)^m - 1}{\beta}\\
     &\leq 2\eta L^2(\mu+1)(\Delta + 2)m.
   \end{aligned}
 \end{equation}
From the definition of $R_t^{s+1}$, we obtain
\begin{equation}
  \begin{aligned}
    R_{t+1}^{s+1} = &\mathbf{E} \Big[ P(x_{t+1}^{s+1}) + c_{t+1}||x_{t+1}^{s+1} - \tilde{x}^s||^2 \Big]\\
    \leq &\mathbf{E} \Big[ P(x_{t+1}^{s+1}) + c_{t+1}(1+\frac{1}{\beta})||x_{t+1}^{s+1} - x_t^{s+1}||^2 + c_{t+1}(1+\beta)||x_t^{s+1}-\tilde{x}^s||^2\Big]\\
    \leq &\mathbf{E}\Big[ P(x_t^{s+1}) + c_t||x_t^{s+1} - \tilde{x}^s||^2 + (c_{t+1}(1+\frac{1}{\beta}) + \frac{L}{2} - \frac{1}{2\eta})||x_{t+1}^{s+1} - x_t^{s+1}||^2 + (L-\frac{1}{2\eta})||\bar{x}_{t+1}^{s+1} - x_t^{s+1}||^2\\
    & + \eta L^2(\mu+1)(\Delta+2)\tau\sum\limits_{d=D(t)}^{t-1} ||x_{d+1}^{s+1} - x_d^{s+1}||^2 \Big].
  \end{aligned}
\end{equation}
Summing over $t = 0$ to $m-1$, we get
\begin{equation}
  \begin{aligned}
    \sum\limits_{t=0}^{m-1}R_{t+1}^{s+1} \leq \sum\limits_{t=0}^{m-1}&R_t^{s+1} + \sum\limits_{t=0}^{m-1}\Big[ c_{t+1}(1+\frac{1}{\beta}) + \frac{L}{2} - \frac{1}{2\eta} + \eta L^2(\mu+1)(\Delta + 2)\tau^2\Big] \mathbf{E}||x_{t+1}^{s+1} - x_{t}^{s+1}||^2 \\
    &+ \sum\limits_{t=0}^{m-1}(L-\frac{1}{2\eta})\mathbf{E}||\bar{x}_{t+1}^{s+1} - x_t^{s+1}||^2.
  \end{aligned}
\end{equation}
The inequality holds because $\sum\limits_{t=0}^{m-1} \sum\limits_{d=D(t)}^{t-1}||x_{d+1}^{s+1} - x_d^{s+1}||^2 \leq \tau\sum\limits_{t=0}^{m-1}||x_{t+1}^{s+1} - x_t^{s+1}||^2$.

\noindent
Now we derive the bound for $\eta$ to make $\Big[c_{t+1}(1+\frac{1}{\beta}) + \frac{L}{2} - \frac{1}{2\eta} + \eta L^2(\mu+1)(\Delta + 2)\tau^2 \Big]\leq 0$.
Since ${c_t}$ is a decreasing sequence, we only need to prove the above inequality for $c_0$.
Let $\eta = \frac{\rho}{L}$, where $\rho < \frac{1}{2}$ is a positive constant.
After calculations, we obtain the following constraint:
\begin{equation}
  \label{eq:tau-inequ1}
     8\rho^2 m^2(\mu+1)(\Delta+2) + 2\rho^2(\mu+1)(\Delta+2)\tau^2 + \rho \leq 1.
\end{equation}
If (\ref{eq:tau-inequ1}) holds, then
\begin{equation}
  \begin{aligned}
    \label{eq:k1}
  \sum\limits_{t=0}^{m-1}  R_{t+1}^{s+1} &\leq \sum\limits_{t=0}^{m-1} R_t^{s+1} + (L-\frac{1}{2\eta})\sum\limits_{t=0}^{m-1}\mathbf{E}||\bar{x}_{t+1}^{s+1} - x_t^{s+1}||^2.
  \end{aligned}
\end{equation}
Because $x_0^{s+1} = \tilde{x}^s$, $x_m^{s+1} = \tilde{x}^{s+1}$ and $c_m = 0$, we have $R(x_0^{s+1}) = P(\tilde{x}^s)$ and $R(x_m^{s+1}) = P(\tilde{x}^{s+1})$.
Summing (\ref{eq:k1}) over $s=0$ to $S-1$, we get
\begin{equation}
  (\frac{1}{2\eta} - L) \sum\limits_{s=0}^{S-1}\sum\limits_{t=0}^{m-1}\mathbf{E}||\bar{x}_{t+1}^{s+1} - x_t^{s+1}||^2 \leq P(x^0) - P(x^*).
\end{equation}
Using the definition of  $G_\eta(x_t^{s+1}) \triangleq \frac{1}{\eta}[x_t^{s+1}-\prox_{\eta h}(x_t^{s+1}-\eta\nabla f(x_t^{s+1}))] =\frac{1}{\eta}(x_t^{s+1} - \bar{x}_{t+1}^{s+1})$, we obtain Theorem $2$.
%
\QEDB

\noindent
\section{Analysis for Sparse-AsyLPG}
%
\begin{lemma}
  \label{lemma:2-spar-variance}
  Define $\varphi_t \triangleq \sum_{i=1}^d p_i$.
  If $\varphi_t \leq \frac{||\alpha_t||_1}{||\alpha_t||_{\infty}}$, then for $\alpha_t = \sum_{i=1}^d \alpha_{t,i} \mathit{e}_i$,
   we have $\mathbf{E}||\beta_t||^2 \geq \frac{1}{\varphi_t}||\alpha_t||^2_1$. 
   The equality holds if and only if $p_i = \frac{|\alpha_{t,i}| \cdot \varphi_t}{||\alpha_t||_1}$.
\end{lemma}
\noindent
\emph{Proof Sketch.}
From the calculation of $\beta_t$, we obtain $\mathbf{E}||\beta_t||^2 = \sum\limits_{i=1}^d\frac{\alpha_{t,i}^2}{p_i}$.
If $\varphi_t \leq \frac{||\alpha_t||_1}{||\alpha_t||_{\infty}}$,
%
then it can be concluded from Lemma $3$ and Theorem $5$ in \cite{wang2018atomo}  that  $\mathbf{E}||\beta_t||^2 \geq \frac{1}{\varphi_t}||\alpha_t||^2_1$, with equality if and only if $p_i = \frac{|\alpha_{t,i}| \cdot \varphi_t}{||\alpha_t||_1}$.
\QEDB

\begin{lemma} 
  \label{lemma:u-f}
  Suppose $\varphi_t \leq \frac{||\alpha_t||_1}{||\alpha_t||_{\infty}}$, Assumptions $1$, $2$, $3$ hold and for each $i \in \{1,...,d\}$, $p_i = \frac{|\alpha_{t,i}|\cdot \varphi_t}{||\alpha_t||_1}$.
  Denote $\Gamma = \frac{d^2}{4\varphi(2^{b-1}-1)^2} + \frac{d}{\varphi} +1$,  where $\varphi = \min_t\{\varphi_t\}$.
  Then, for the gradient $u_{t}^{s+1}$ in Algorithm $2$, we have
  \begin{equation}
    \begin{aligned}
      \label{eq:saprse-u-var}
      \mathbf{E}||u_t^{s+1} - \nabla f(x_t^{s+1})||^2 \leq
      2L^2(\mu+1)\Gamma \mathbf{E} \Big[ ||x_{D(t)}^{s+1} - x_t^{s+1}||^2 + ||x_t^{s+1} - \tilde{x}^s||^2 \Big].
\end{aligned}
\end{equation}
\end{lemma}
\noindent
\emph{Proof.}
\begin{equation}
  \begin{aligned}
    & \mathbf{E}||u_t^{s+1} - \nabla f(x_t^{s+1})||^2 \\
    & = \mathbf{E}||Q_{(\delta_{\beta_t},b)}(\beta_t) + \nabla f(\tilde{x}^s) - \nabla f(x_t^{s+1})||^2\\
    & = \mathbf{E}||Q_{(\delta_{\beta_t},b)}(\beta_t) - \beta_t||^2 + \mathbf{E}||\beta_t + \nabla f(\tilde{x}^s) - \nabla f(x_t^{s+1})||^2 \\
    & \leq \frac{d}{4(2^{b-1}-1)^2}\mathbf{E}||\beta_t||^2 + \mathbf{E}||\beta_t - \alpha_t||^2 + \mathbf{E}||\alpha_t + \nabla f(\tilde{x}^s) - \nabla f(x_t^{s+1})||^2\\
    &= \Big[ \frac{d}{4(2^{b-1}-1)^2} + 1\Big]\mathbf{E}||\beta_t||^2 - \mathbf{E}||\alpha_t||^2 +  \mathbf{E}||\alpha_t + \nabla f(\tilde{x}^s) - \nabla f(x_t^{s+1})||^2
    \end{aligned}
\end{equation}
where the second equality holds because $Q_{(\delta_{\beta_t},b)}(\beta_t)$ is an unbiased quantization.
The first inequality uses Lemma~\ref{lemma:unb-quan-var} and $\mathbf{E}[\beta_t] = \alpha_t$.
According to $T_2$ we get
\begin{equation}
  \begin{aligned}
& \mathbf{E}||u_t^{s+1} - \nabla f(x_t^{s+1})||^2 \\
&\leq \Big[ \frac{d}{4(2^{b-1}-1)^2} + 1\Big] \mathbf{E}||\beta_t||^2 - \mathbf{E}||\alpha_t||^2 + 4L^2(\mu+1)\mathbf{E}||x_{D(t)}^{s+1} - x_t^{s+1}||^2 + 4L^2(\mu+1)\mathbf{E}||x_t^{s+1} - \tilde{x}^s||^2.
  \end{aligned}
\end{equation}
From Lemma \ref{lemma:2-spar-variance}, we obtain
\begin{equation}
  \begin{aligned}
    \label{eq:u-f-2}
  & \mathbf{E}||u_t^{s+1} - \nabla f(x_t^{s+1})||^2 \\
&\leq \Big[ \frac{d^2}{4\varphi(2^{b-1}-1)^2} + \frac{d}{\varphi} - 1\Big] \mathbf{E}||\alpha_t||^2 + 4L^2(\mu+1)\mathbf{E}||x_{D(t)}^{s+1} - x_t^{s+1}||^2 + 4L^2(\mu+1)\mathbf{E}||x_t^{s+1} - \tilde{x}^s||^2\\
&\leq 2L^2(\mu+1)\Big[ \frac{d^2}{4\varphi(2^{b-1}-1)^2} + \frac{d}{\varphi} +1 \Big] \mathbf{E}||x_{D(t)}^{s+1} - x_t^{s+1}||^2 + 2L^2(\mu+1)\Big[ \frac{d^2}{4\varphi(2^{b-1}-1)^2} + \frac{d}{\varphi} +1 \Big] \mathbf{E} ||x_t^{s+1} - \tilde{x}^s||^2,
  \end{aligned}
\end{equation}
where the first inequality uses $||x||_1 \leq \sqrt{d}||x||_2$ for $x \in \mathbb{R}^d$ and $\mathbf{E}||\beta_t||^2 = \frac{1}{\varphi_t}||\alpha_t||^2_1$ when $p_i = \frac{|\alpha_{t,i}|\cdot \varphi_t}{||\alpha_t||_1}$.
The final inequality comes from the upper bound of $T_1$.
\QEDB
\\
\\
\noindent
\emph{\textbf{Proof sketch of Theorem $3$.}}
Substituting (\ref{eq:saprse-u-var}) in (\ref{eq:P-in}) and following the proof of Theorem $2$, we obtain a convergence rate of
\begin{equation}
  \mathbf{E}||G_\eta(x_{out})||^2 \leq \frac{2L[P(x^0) - P(x^*)]}{\rho(1-2\rho)T},
\end{equation}
if  $8\rho^2 m^2(\mu+1)\Gamma + 2\rho^2(\mu+1)\tau^2\Gamma + \rho \leq 1$, where $\Gamma = \frac{d^2}{4\varphi(2^{b-1}-1)^2} + \frac{d}{\varphi} +1$.
\QEDB

\section{Proof of Theorem $4$}
The following lemma is a widely used technical result in composite optimization, which is called \emph{3-Point-Property}. Lemma $1$ in \cite{Lan2012An} provides its detailed proofs and extensions.
\begin{lemma}\label{lemma:3-point property}
   If $y_{t+1}^s$ is the optimal solution of
  \begin{equation}
    \min_{y \in \chi} \phi(y) + \frac{1}{2\eta_s}||y-y_{t}^s||^2,
  \end{equation}
  where function $\phi(y)$ is convex over a convex set $\chi$. Then for any $y \in \chi$, we have \cite{Lan2012An}
  \begin{equation}
    \label{eq:3-point}
    \phi(y) + \frac{1}{2\eta_s}||y-y_{t}^s||^2 \geq \phi(y_{t+1}^s) + \frac{1}{2\eta_s}||y_{t+1}^s - y_{t}^s||^2 + \frac{1}{2\eta_s}||y-y_{t+1}^s||^2.
  \end{equation}
\end{lemma}

\noindent
\emph{\textbf{Proof of Theorem $4$}.}
From the update rule of $y_{t+1}^s$, we know that
\begin{equation}
  y_{t+1}^s = \argmin_{y} h(y)+\langle u_t^s, y-y_t^s\rangle + \frac{1}{2\eta_s}||y-y_t^s||^2.
\end{equation}
Applying Lemma \ref{lemma:3-point property} with $\phi(y) = h(y) + \langle u_t^s, y-y_t^s\rangle$ and $y=x^*$ in (\ref{eq:3-point}), we obtain
\begin{equation}
  \label{eq:asy-msvrg-3-point}
  h(y_{t+1}^s) + \langle u_t^s, y_{t+1}^s - y_t^s\rangle \leq h(x^*) + \langle u_t^s, x^* - y_t^s\rangle + \frac{1}{2\eta_s}||x^* - y_t^s||^2 - \frac{1}{2\eta_s}||x^* - y_{t+1}^s||^2 - \frac{1}{2\eta_s}||y_{t+1}^s - y_t^s||^2.
\end{equation}
Since $f(x)$ is Lipschitz smooth, we have
\begin{equation}
  \begin{aligned}
    \label{eq:asy-msvrg-f-lip}
  \mathbf{E} f(x_{t+1}^s) \leq & \mathbf{E} \Big( f(x_t^s) + \langle \nabla f(x_t^s), x_{t+1}^s - x_t^s\rangle + \frac{L}{2}||x_{t+1}^s - x_t^s||^2 \Big)\\
  = & \mathbf{E} \Big(f(x_t^s) + \theta_s \langle u_t^s, y_{t+1}^s - y_t^s\rangle + \theta_s\langle \nabla f(x_t^s) - u_t^s, y_{t+1}^s - y_t^s\rangle + \frac{L}{2}||x_{t+1}^s - x_t^s||^2 \Big),
\end{aligned}
\end{equation}
where the first equality uses $x_{t+1}^s - x_t^s = \theta_s(y_{t+1}^s - y_t^s)$. Therefore,
\begin{equation}
  \begin{aligned}
    \label{eq:asy-m-P-1}
    \mathbf{E} P(x_{t+1}^s) = &\mathbf{E} \Big[ f(x_{t+1}^s) + h(x_{t+1}^s)\Big]\\
    \leq & \mathbf{E} \Big[ f(x_t^s) + \theta_s \langle u_t^s, y_{t+1}^s - y_t^s\rangle + \theta_s\langle \nabla f(x_t^s) - u_t^s, y_{t+1}^s - y_t^s\rangle + \frac{L}{2}||x_{t+1}^s - x_t^s||^2 + h(x_{t+1}^s)  \Big]\\
    \leq & \mathbf{E} \Big[ f(x_t^s) + \theta_s \langle u_t^s, y_{t+1}^s - y_t^s\rangle + \theta_s\langle \nabla f(x_t^s) - u_t^s, y_{t+1}^s - y_t^s\rangle + \frac{L}{2}||x_{t+1}^s - x_t^s||^2 + (1-\theta_s)h(\tilde{x}^{s-1}) + \theta_sh(y_{t+1}^s)  \Big]\\
    \leq & \mathbf{E} \Big[ \theta_s h(x^*) + \underbrace{\theta_s\langle u_t^s, x^* - y_t^s\rangle}_{T_3} + \frac{\theta_s}{2\eta_s}||x^* - y_t^s||^2 - \frac{\theta_s}{2\eta_s}||x^* - y_{t+1}^s||^2 - \frac{\theta_s}{2\eta_s}||y_{t+1}^s - y_t^s||^2 \\
    & + f(x_t^s) + \underbrace{\theta_s\langle \nabla f(x_t^s) - u_t^s, y_{t+1}^s - y_t^s\rangle}_{T_4} + \frac{L}{2}||x_{t+1}^s - x_t^s||^2 + (1-\theta_s)h(\tilde{x}^{s-1}) \Big],
  \end{aligned}
\end{equation}
where the first inequality uses (\ref{eq:asy-msvrg-f-lip}), and the second inequality follows from $x_{t+1}^s = \theta_s y_{t+1}^s + (1-\theta_s)\tilde{x}^{s-1}$ and the convexity of $h(x)$. We apply (\ref{eq:asy-msvrg-3-point}) in the third inequality.
$T_3$ can be bounded as follows.
\begin{equation}
  \begin{aligned}
  \mathbf{E} T_3 &= \theta_s \mathbf{E} \langle u_t^s, x^* - y_t^s\rangle\\
  &=\mathbf{E} \langle u_t^s, \theta_s x^* + (1-\theta_s)\tilde{x}^{s-1} - x_t^s\rangle\\
  &=\mathbf{E} \langle u_t^s, \theta_s x^* + (1-\theta_s)\tilde{x}^{s-1} - Q_{(\delta_x, b_x)}(x_{D(t)}^s)\rangle + \mathbf{E}\langle u_t^s, Q_{(\delta_x, b_x)}(x_{D(t)}^s) - x_t^s\rangle\\
  &=\mathbf{E} \langle \nabla f_a(Q_{(\delta_x, b_x)}(x_{D(t)}^s)), \theta_s x^* + (1-\theta_s)\tilde{x}^{s-1} - Q_{(\delta_x, b_x)}(x_{D(t)}^s)\rangle + \mathbf{E}\langle \nabla f_a(Q_{(\delta_x, b_x)}(x_{D(t)}^s), Q_{(\delta_x, b_x)}(x_{D(t)}^s) - x_t^s\rangle\\
  &\leq \mathbf{E} \Big[ f_a(\theta_s x^* + (1-\theta_s)\tilde{x}^{s-1}) - f_a(Q_{(\delta_x, b_x)}(x_{D(t)}^s)) + f_a(Q_{(\delta_x, b_x)}(x_{D(t)}^s)) - f_a(x_t^s) + \frac{L}{2}||Q_{(\delta_x, b_x)}(x_{D(t)}^s) - x_t^s||^2 \Big]\\
  &\leq \mathbf{E} \Big[ \theta_s f(x^*) + (1-\theta_s)f(\tilde{x}^{s-1}) - f(x_t^s) + \frac{L}{2}||Q_{(\delta_x, b_x)}(x_{D(t)}^s) - x_t^s||^2\Big],
\end{aligned}
\end{equation}
where the convexity and Lipschitz smoothness of $f_a(x)$ are adopted in the first inequality.
Next we derive the bound of $\mathbf{E}||Q_{(\delta_x, b_x)}(x_{D(t)}^s) - x_t^s||^2$ as follows.
\begin{equation}
  \begin{aligned}
    \mathbf{E}||Q_{(\delta_x, b_x)}(x_{D(t)}^s) - x_t^s||^2
    &= \mathbf{E}||Q_{(\delta_x, b_x)}(x_{D(t)}^s) - x_{D(t)}^s||^2 + \mathbf{E}||x_{D(t)}^s - x_t^s||^2\\
    &\leq \theta_s\mu\mathbf{E}||x_{D(t)}^s - \tilde{x}^{s-1}||^2 + \mathbf{E}||x_{D(t)}^s - x_t^s||^2\\
    &\leq (1+2\theta_s\mu)\mathbf{E}||x_{D(t)}^s - x_t^s||^2 + 2\theta_s^3 \mu\mathbf{E} ||y_t^s - \tilde{x}^{s-1}||^2\\
    &\leq (1+2\theta_s\mu)\mathbf{E}||x_{D(t)}^s - x_t^s||^2 + 2\theta_s^3\mu D.
  \end{aligned}
\end{equation}
where the first equality holds because $Q_{(\delta_x, b_x)}(x_{D(t)}^s)$ is an unbiased quantization and
the first inequality comes from Step $8$ in Algorithm $3$.
The second inequality holds because $x_t^s - \tilde{x}^{s-1} = \theta_s(y_t^s - \tilde{x}^{s-1})$.
Therefore,
\begin{equation}
  \label{eq:T_1}
  \mathbf{E}T_3 \leq \mathbf{E} \Big[  \theta_s f(x^*) + (1-\theta_s)f(\tilde{x}^{s-1}) - f(x_t^s) + \frac{(1+2\theta_s\mu)L}{2}||x_{D(t)}^s - x_t^s||^2 + \theta_s^3\mu LD\Big].
\end{equation}
Now we bound $T_4$.
Define $v_t^s  \triangleq \nabla f_a(x_t^s) - \nabla f_a(\tilde{x}^{s-1}) + \nabla f(\tilde{x}^{s-1})$.
\begin{equation}
\mathbf{E}T_4 = \theta_s\mathbf{E}\langle \nabla f(x_t^s) - u_t^s, y_{t+1}^s - y_t^s\rangle = \underbrace{\theta_s\mathbf{E}\langle \nabla f(x_t^s) - v_t^s, y_{t+1}^s - y_t^s\rangle}_{T_5} + \underbrace{\theta_s\mathbf{E}\langle v_t^s - u_t^s, y_{t+1}^s - y_t^s\rangle}_{T_6}.
\end{equation}
\begin{equation}
  \begin{aligned}
    \label{eq:asyn-mom-T_3}
   T_5 =& \theta_s\mathbf{E}\langle \nabla f(x_t^s) - v_t^s, y_{t+1}^s - y_t^s\rangle\\
  \leq& \frac{\theta_s}{2\tau L}\mathbf{E}||\nabla f(x_t^s) - v_t^s||^2 + \frac{\tau L \theta_s}{2}\mathbf{E}||y_{t+1}^s - y_t^s||^2\\
  \leq& \frac{\theta_s}{2\tau L}\mathbf{E}||\nabla f_{a}(x_t^s) - \nabla f_{a}(\tilde{x}^{s-1})||^2 + \frac{\tau L \theta_s}{2}\mathbf{E}||y_{t+1}^s - y_t^s||^2\\
  \leq& \frac{\theta_sL^2}{2\tau L}\mathbf{E}||x_t^s - \tilde{x}^{s-1}||^2 + \frac{\tau L \theta_s}{2}\mathbf{E}||y_{t+1}^s - y_t^s||^2\\
  =& \frac{\theta_s^3 L^2}{2\tau L}\mathbf{E}||y_t^s - \tilde{x}^{s-1}||^2 + \frac{\tau L \theta_s}{2}\mathbf{E}||y_{t+1}^s - y_t^s||^2\\
  \leq& \frac{\theta_s^3 LD}{2\tau} + \frac{\tau L \theta_s}{2}\mathbf{E}||y_{t+1}^s - y_t^s||^2,
\end{aligned}
\end{equation}
where in the first inequality we use Young's inequality.
The second equality follows from $x_t^s - \tilde{x}^{s-1} = \theta_s (y_t^s - \tilde{x}^{s-1})$.
Moreover,
\begin{equation}
  \begin{aligned}
    \label{eq:asyn-mom-T_4}
    T_6 &=  \theta_s\mathbf{E}\langle v_t^s - u_t^s, y_{t+1}^s - y_t^s\rangle\\
    &\leq \frac{\theta_s}{2\tau L}\mathbf{E}||v_t^s - u_t^s||^2 + \frac{\tau L\theta_s }{2}\mathbf{E}||y_{t+1}^s - y_t^s||^2.
  \end{aligned}
\end{equation}
From the definition of $u_t^s$ and $v_t^s$, we have
\begin{equation}
  \begin{aligned}
    &\mathbf{E}||v_t^s - u_t^s||^2 \\
    &= \mathbf{E} ||Q_{(\delta_{\alpha_t},b)}(\alpha_t)- \nabla f_a(x_t^s) + \nabla f_a(\tilde{x}^{s-1})||^2\\
    &= \mathbf{E}||Q_{(\delta_{\alpha_t},b)} \Big( \nabla f_a(Q_{(\delta_x, b_x)}(x_{D(t)}^s)) - \nabla f_a(\tilde{x}^{s-1})\Big) - \nabla f_a(Q_{(\delta_x, b_x)}(x_{D(t)}^s)) + \nabla f_a(\tilde{x}^{s-1})||^2 \\
    & \quad + \mathbf{E} || \nabla f_a(Q_{(\delta_x, b_x)}(x_{D(t)}^s)) - \nabla f_a(x_t^s)||^2\\
    &\leq \frac{d}{4(2^{b-1}-1)^2}\mathbf{E}||\nabla f_a(Q_{(\delta_x, b_x)}(x_{D(t)}^s)) - \nabla f_a(\tilde{x}^{s-1})||^2 + \mathbf{E} || \nabla f_a(Q_{(\delta_x, b_x)}(x_{D(t)}^s)) - \nabla f_a(x_t^s)||^2\\
    &\leq \frac{dL^2}{4(2^{b-1}-1)^2}\mathbf{E}||Q_{(\delta_x, b_x)}(x_{D(t)}^s) - \tilde{x}^{s-1}||^2 + L^2\mathbf{E}||Q_{(\delta_x, b_x)}(x_{D(t)}^s) - x_t^s||^2\\
    &\leq \Big[ \frac{dL^2}{2(2^{b-1}-1)^2} + L^2\Big] \mathbf{E}||Q_{(\delta_x, b_x)}(x_{D(t)}^s)  - x_t^s||^2 + \frac{dL^2}{2(2^{b-1}-1)^2} \mathbf{E}||x_t^s - \tilde{x}^{s-1}||^2\\
    &\leq \Big[ \frac{dL^2}{(2^{b-1}-1)^2} + 2L^2\Big] \mathbf{E}||Q_{(\delta_x, b_x)}(x_{D(t)}^s)  - x_{D(t)}^s||^2 + \Big[ \frac{dL^2}{(2^{b-1}-1)^2} + 2L^2\Big] \mathbf{E}||x_{D(t)}^s - x_t^s||^2 \\
     & \quad \quad \quad + \frac{dL^2\theta_s^2}{2(2^{b-1}-1)^2}\mathbf{E}||y_t^{s} - \tilde{x}^{s-1}||^2\\
    &\leq \Big[ \frac{dL^2}{(2^{b-1}-1)^2} + 2L^2\Big]\theta_s \mu \mathbf{E}||x_{D(t)}^s - \tilde{x}^{s-1}||^2 + \Big[ \frac{dL^2}{(2^{b-1}-1)^2} + 2L^2\Big] \mathbf{E}||x_{D(t)}^s - x_t^s||^2 + \frac{dL^2\theta_s^2D}{2(2^{b-1}-1)^2}\\
    &\leq \Big[ \frac{dL^2}{(2^{b-1}-1)^2} + 2L^2\Big](1+2\theta_s\mu)\mathbf{E}||x_{D(t)}^s - x_t^s||^2 + \Big[ \frac{2dL^2}{(2^{b-1}-1)^2} + 4L^2\Big]\theta_s\mu \mathbf{E}||x_t^s - \tilde{x}^{s-1}||^2 + \frac{dL^2\theta_s^2D}{2(2^{b-1}-1)^2}\\
    &\leq (1+2\theta_s\mu)L^2\Big[ \frac{d}{(2^{b-1}-1)^2} + 2\Big] \mathbf{E}||x_{D(t)}^s - x_t^s||^2 + \theta_s^3L^2\mu\Big[ \frac{2d}{(2^{b-1}-1)^2} + 4\Big]D +  \frac{dL^2\theta_s^2D}{2(2^{b-1}-1)^2},
  \end{aligned}
\end{equation}
where the second equality holds because $Q_{(\delta_{\alpha_t},b)}(\alpha_t)$ is an unbiased quantization,
and the first inequality uses Lemma~\ref{lemma:unb-quan-var}.
In the fourth and final inequality, we adopt $x_t^s - \tilde{x}^{s-1} = \theta_s(y_t^s - \tilde{x}^{s-1})$.
The fifth inequality follows from Step $8$ of Algorithm $3$.
Putting them together, we obtain
\begin{equation}
  \begin{aligned}
  \label{eq:T_2}
  \mathbf{E} T_4 \leq &\tau L \theta_s\mathbf{E}||y_{t+1}^s - y_t^s||^2 + \frac{L\theta_s(1+2\theta_s\mu)}{2\tau}\Big[ \frac{d}{(2^{b-1}-1)^2} + 2\Big] \mathbf{E}||x_{D(t)}^s - x_t^s||^2 + \frac{LD\theta_s^3}{2\tau}\Big[\frac{d}{2(2^{b-1}-1)^2} + 1\Big] \\
  &+ \frac{\theta_s^4LD\mu}{2\tau}\Big[ \frac{2d}{(2^{b-1}-1)^2} + 4\Big].
\end{aligned}
\end{equation}
Substituting (\ref{eq:T_2}) and (\ref{eq:T_1}) in (\ref{eq:asy-m-P-1}), we get
\begin{equation}
  \begin{aligned}
    \mathbf{E} P(x_{t+1}^s) \leq &\mathbf{E} \Big[  (1-\theta_s)P(\tilde{x}^{s-1}) + \theta_s P(x^*) + \frac{\theta_s}{2\eta_s}(||x^* - y_t^s||^2 - ||x^* - y_{t+1}^s||^2)\\
    &+ \theta_s^3\mu LD +  \frac{LD\theta_s^3}{2\tau}\Big[\frac{d}{2(2^{b-1}-1)^2} + 1\Big] + \frac{\theta_s^4LD\mu}{2\tau}\Big[ \frac{2d}{(2^{b-1}-1)^2} + 4\Big]\\
    & + \underbrace{\frac{(1+2\theta_s\mu)L}{2}||x_{D(t)}^s - x_t^s||^2 + \tau L \theta_s||y_{t+1}^s - y_t^s||^2 + \frac{L\theta_s(1+2\theta_s\mu)}{2\tau}\Big[ \frac{d}{(2^{b-1}-1)^2} + 2\Big]||x_{D(t)}^s - x_t^s||^2  }_{T_7} \\
    & + \underbrace{\frac{L}{2}||x_{t+1}^s - x_t^s||^2 - \frac{\theta_s}{2\eta_s}||y_{t+1}^s - y_t^s||^2}_{T_8}\Big].
\end{aligned}
\end{equation}
Let $\eta_s\theta_s = \frac{1}{\sigma L}$ where $\sigma > 1$, since $\sum\limits_{t=0}^{m-1}\sum\limits_{d=D(t)}^{t-1}||x_{d+1}^s - x_d^s||^2 \leq \tau \sum\limits_{t=0}^{m-1}||x_{t+1}^s - x_t^s||^2$, it can be verified that
\begin{equation}
  \sum\limits_{t=0}^{m-1}(T_7+T_8) \leq \xi\sum\limits_{t=0}^{m-1}||y_{t+1}^s - y_t^s||^2,
\end{equation}
where $\xi = \tau^2\theta_s^2\Big[\frac{(1+2\theta_s\mu)L}{2} + \frac{(1+2\theta_s\mu)\theta_s L}{2\tau}(\frac{d}{(2^{b-1}-1)^2} + 2)\Big] + \tau L \theta_s + \frac{L\theta_s^2}{2} - \frac{\sigma L \theta_s^2}{2}$.
Denote $\Delta = \frac{d}{(2^{b-1}-1)^2} + 2$,
if $\tau \leq \frac{ \sqrt{ \big(\frac{2}{(1+2\theta_s\mu)\theta_s} + \theta_s\Delta \big)^2 + \frac{4(\sigma -1)}{(1+2\theta_s\mu)}} - (\frac{2}{(1+2\theta_s\mu)\theta_s} + \theta_s\Delta)}{2}$, then $\xi \leq 0$.
Suppose the above constraint holds,  we have
\begin{equation}
  \sum\limits_{t=0}^{m-1}\mathbf{E} P(x_{t+1}^s) \leq \sum\limits_{t=0}^{m-1}\mathbf{E} \Big[ (1-\theta_s)P(\tilde{x}^{s-1}) + \theta_s P(x^*) + \frac{\theta_s}{2\eta_s}(||x^* - y_t^s||^2 - ||x^* - y_{t+1}^s||^2) + \theta_s^3\mu LD + \frac{\theta_s^3LD\Delta}{4\tau} + \frac{\theta_s^4LD\Delta\mu}{\tau}\Big].
\end{equation}
Using $\tilde{x}^s = \frac{1}{m}\sum\limits_{t=0}^{m-1}x_{t+1}^s$, we obtain
\begin{equation}
  \label{eq:asy-mom-F}
  \mathbf{E} \Big[ P(\tilde{x}^s) - P(x^*)\Big] \leq \mathbf{E} \Big[ (1-\theta_s)\big(P(\tilde{x}^{s-1}) - P(x^*)\big) + \frac{\sigma L \theta_s^2}{2m}(||y_{0}^{s} - x^*||^2 - ||y_m^s - x^*||^2) + \theta_s^3\mu LD + \frac{\theta_s^3LD\Delta}{4\tau}  + \frac{\theta_s^4LD\Delta\mu}{\tau}\Big].
\end{equation}
Dividing both sides of (\ref{eq:asy-mom-F}) by $\theta_s^2$, summing over $s=1$ to $S$, and using the definition $y_0^s = y_{m}^{s-1}$ and that $\frac{1-\theta_s}{\theta_s^2}\leq \frac{1}{\theta_{s-1}^2}$ when $\theta_s=\frac{2}{s+2}$,  we have
\begin{equation}
  \label{eq:asy-msvrg-final}
  \mathbf{E}\Big[ P(\tilde{x}^S) - P(x^*)\Big] \leq \frac{4\Big[P(\tilde{x}^0) - P(x^*)\Big]}{(S+2)^2} + \frac{2\sigma L||\tilde{x}^{0}-x^*||^2}{m(S+2)^2} + \frac{8\big(\Delta(1+\mu)/\tau + \mu \big)LD\log(S+2)}{(S+2)^2}.
\end{equation}
\QEDB



\end{document}